\documentclass[12pt,leqno]{article}
\usepackage{amsmath, amsthm, amsfonts, amssymb, color}
\usepackage{mathrsfs}
\newtheorem{thm}{Theorem}[section]

\newtheorem{lem}[thm]{Lemma}
\newtheorem{prp}[thm]{Proposition}

\newtheorem{exa}[thm]{Example}
\theoremstyle{definition}

\newcommand{\scr}[1]{\mathscr #1}
\definecolor{wco}{rgb}{0.5,0.2,0.3}

\numberwithin{equation}{section} \theoremstyle{remark}
\newtheorem{rem}{Remark}[section]

\newcommand{\ua}{\uparrow}

\title{{\bf  Stochastic Generalized Porous Media and Fast Diffusion
Equations} \footnote{Supported in part by the DFG through the
Forschergruppe ``Spectral Analysis, Asymptotic Distributions and
Stochastic Dynamics'',  the SFB 701, the BiBoS Research Centre,
NNSFC(10121101) and RFDP(20040027009).} }
\author{{\bf Jiagang Ren}\\
\footnotesize{School of Mathematics and Computational Science,
Zhongshan University} \\ \footnotesize{Guangzhou, Guangdong
510275, China}\\
\\ {\bf  Michael R\"ockner}\\
\footnotesize{Departments of Mathematics and Statistics, Purdue University,
West Lafayette, IN, 47906, USA}\\
{\bf Feng-Yu Wang\footnote{The corresponding author: wangfy@bnu.edu.cn}}\\
\footnotesize{School of Mathematical Sciences, Beijing Normal
University, Beijing 100875, China}\\
}

\begin{document}
\maketitle
\begin{abstract} We present a generalization of  Krylov-Rozovskii's
result on the existence and uniqueness of solutions to  monotone
stochastic differential equations. As an application, the
stochastic generalized porous media and fast diffusion  equations
are studied for $\sigma$-finite  reference measures, where the
drift term is given by a negative definite operator acting on a
time-dependent function, which belongs to a large class of
functions comparable with the so-called $N$-functions in the
theory of Orlicz spaces.
\end{abstract} \noindent
 \vskip 2cm

\def\R{\mathbb R}  \def\ff{\frac} \def\ss{\sqrt} \def\BB{\mathbb
B}
\def\N{\mathbb N} \def\kk{\kappa} \def\m{{\bf m}}
\def\dd{\delta} \def\DD{\Delta} \def\vv{\varepsilon} \def\rr{\rho}
\def\<{\langle} \def\>{\rangle} \def\GG{\Gamma} \def\gg{\gamma}
  \def\nn{\nabla} \def\pp{\partial} \def\tt{\tilde}
\def\d{\text{\rm{d}}} \def\bb{\beta} \def\aa{\alpha} \def\D{\scr D}
\def\E{\mathbb E} \def\si{\sigma} \def\ess{\text{\rm{ess}}}
\def\beg{\begin} \def\beq{\begin{equation}}  \def\F{\scr F}
\def\Ric{\text{\rm{Ric}}} \def\Hess{\text{\rm{Hess}}}\def\B{\mathbb B}
\def\e{\text{\rm{e}}} \def\ua{\underline a} \def\OO{\Omega} \def\b{\mathbf b}
\def\oo{\omega}     \def\tt{\tilde} \def\Ric{\text{\rm{Ric}}}
\def\cut{\text{\rm{cut}}} \def\P{\mathbb P} \def\ifn{I_n(f^{\bigotimes n})}
\def\fff{f(x_1)\dots f(x_n)} \def\ifm{I_m(g^{\bigotimes m})} \def\ee{\varepsilon}
\def\pm{\pi_{{\bf m}}}   \def\p{\mathbf{p}}   \def\ml{\mathbf{L}}
 \def\C{\scr C}      \def\aaa{\mathbf{r}}     \def\r{r}
\def\gap{\text{\rm{gap}}} \def\prr{\pi_{{\bf m},\varrho}}  \def\r{\mathbf r}
\def\Z{\mathbb Z} \def\vrr{\varrho} \def\ll{\lambda}
\def\L{\scr L}


\newcommand*{\argdot}{\,\cdot\,}
\newcommand*{\eps}{\varepsilon}


\newcommand*{\norm}[1]{\lVert#1\rVert}
\newcommand*{\lrnorm}[1]{\left\lVert#1\right\rVert}
\newcommand*{\bignorm}[1]{\bigl\lVert#1\bigr\rVert}
\newcommand*{\Bignorm}[1]{\Bigl\lVert#1\Bigr\rVert}
\newcommand*{\biggnorm}[1]{\biggl\lVert#1\biggr\rVert}
\newcommand*{\Biggnorm}[1]{\Biggl\lVert#1\Biggr\rVert}

\newcommand*{\abs}[1]{\lvert#1\rvert}
\newcommand*{\lrabs}[1]{\left\lvert#1\right\rvert}
\newcommand*{\bigabs}[1]{\bigl\lvert#1\bigr\rvert}
\newcommand*{\Bigabs}[1]{\Bigl\lvert#1\Bigr\rvert}
\newcommand*{\biggabs}[1]{\biggl\lvert#1\biggr\rvert}
\newcommand*{\Biggabs}[1]{\Biggl\lvert#1\Biggr\rvert}

\newcommand*{\bracket}[1]{\langle#1\rangle}
\newcommand*{\lrbracket}[1]{\left\langle#1\right\rangle}
\newcommand*{\bigbracket}[1]{\bigl\langle#1\bigr\rangle}
\newcommand*{\Bigbracket}[1]{\Bigl\langle#1\Bigr\rangle}
\newcommand*{\biggbracket}[1]{\biggl\langle#1\biggr\rangle}
\newcommand*{\Biggbracket}[1]{\Biggl\langle#1\Biggr\rangle}

\newcommand*{\bottom}[1]{\lfloor#1\rfloor}
\newcommand*{\lrbottom}[1]{\left\lfloor#1\right\rfloor}
\newcommand*{\bigbottom}[1]{\bigl\lfloor#1\bigr\rfloor}
\newcommand*{\Bigbottom}[1]{\Bigl\lfloor#1\Bigr\rfloor}
\newcommand*{\biggbottom}[1]{\biggl\lfloor#1\biggr\rfloor}
\newcommand*{\Biggbottom}[1]{\Biggl\lfloor#1\Biggr\rfloor}

\newcommand*{\pairing}[3]{\sideset{_{#1}^{}}{_{#3}^{}}{\mathop{\bracket{{#2}}}}}
\newcommand*{\bigpairing}[3]{\sideset{_{#1}^{}}{_{#3}^{}}{\mathop{\bigbracket{{#2}}}}}

\newcommand*{\rs}[1]{{\vert}_{#1}}
\newcommand*{\bigrs}[1]{{\big\vert}_{#1}}
\newcommand*{\Bigrs}[1]{{\Big\vert}_{#1}}
\newcommand*{\biggrs}[1]{{\bigg\vert}_{#1}}
\newcommand*{\Biggrs}[1]{{\Bigg\vert}_{#1}}

\section{Introduction}

The main aim of this paper is to solve stochastic partial
differential equations (SPDE) of ``porous media'' type, i.e.
\begin{equation}\label{eq-1.1}
  \mathrm dX_t
  = \bigl( L\Psi(t,X_t) + \Phi(t,X_t) \bigr)\;\mathrm dt
  + B(t,X_t)\;\mathrm dW_t\;,
\end{equation}
where $L$ is a partial (or pseudo) differential operator of order
(less than or equal to) two, so e.g. $L=\Delta$ (or $L=
-(-\Delta)^\alpha$, $\alpha\in(0,1]$) on $\R^d$ or an open subset
thereof (see Example~\ref{exa-3.3} below). The maps
$\Psi,\Phi:[0,T]\times\R\to\R$ (possibly also random) are to fulfill
certain monotonicity conditions (cf.\ \textbf{(A1)}, \textbf{(A2)}
in Section~\ref{sec:stoch-equations} below). $B$ is a
Hilbert-Schmidt operator valued Lipschitz map and $W_t$ a Brownian
motion on a suitable Hilbert space of generalized functions (see
below for details). A typical example for $\Psi$ is
\[
  \Psi(s)
  = \sum_{i=1}^m \alpha_i\,\operatorname{sign}(s)\,\abs s^{r_i}\;,
  \quad s\in \R\;,
\]
with pairwise distinct $r_1,\dots,r_m>0$ and
$\alpha_i\in(0,\infty)$. For $r_i>1$ this corresponds to the
classical porous medium equation and for $r_i\in(0,1)$ to the
so-called fast diffusion equation. Their behaviour is quite
different. E.g.\ in the deterministic case, when $L$ is the
Dirichlet Laplacian on an open bounded domain in $\R^d$, it is
well-known that if $r_i>1$ the solution decays algebraically fast in
$t$, and if $r_i\in(0,1)$, it decays to zero in finite time. We
refer to \cite{Aronson,AP,E} and the references therein, also for
historical remarks.

In recent years, the stochastic version of the porous medium
equation has been studied intensively, see \cite{DRRW} for the
existence, uniqueness   and long-time behavior of some stochastic
generalized porous media equations with finite reference measures,
\cite{Kim} for the stochastic porous media equation on $\R^d$
where the reference (Lebesgue) measure is infinite and \cite{BBDR,
P1} for the analysis of the corresponding Kolmogorov equations.
See also \cite{RWW} for large deviations for a class of
generalized porous media equations.

Our analysis of \eqref{eq-1.1} is based on the so-called variational
approach and requires monotonicity assumptions on the coefficients.
More precisely, we extend classical (impressive) work by
N.V.~Krylov and B.~Rozovskii \cite{KR} (see in particular
Theorems II.2.1 and II.2.2 there and also the pioneering work by
E.~Pardoux \cite{P1,P2}) to make it applicable to SPDE \eqref{eq-1.1}
for very general non-linearities $\Psi$ and $\Phi$.
The corresponding abstract result is proved in Section~\ref{sec:mon-sdes}
below (cf.\ in particular Theorem~\ref{T2.1}) where we also
describe the new framework.
In addition, we include a detailed proof of the crucial It\^o formula,
proved in \cite[Theorem~I.3.1]{KR}, adapted to our more general
situation, in the Appendix.

One of the main points in applying the variational approach is to
find a suitable Gelfand triple
\[
  V\subset H\subset V^*\;,
\]
where $V$ is a separable reflexive Banach space and $V^*$ its dual,
and $H$ is a Hilbert space. In the case of \eqref{eq-1.1} it turns
out that the appropriate Hilbert space is just the Green space of
the operator $L$, i.e.\ the dual of the zero-order Dirichlet space
determined by $L$ (see \eqref{eq-3.0} below and, in particular,
Proposition~\ref{prp-3.1}). For this to be well-defined we need that
the semigroup generated by $L$ is transient. If $L=\Delta$ on a
bounded open subset $\Lambda$ of $\R^d$ with Dirichlet boundary
conditions, then $H=H^{-1}(\Lambda)$ (i.e.\ the dual space of the
Sobolev space $H_0^{1,2}(\Lambda)$). But also the case
$L=-(-\Delta)^\alpha$, $\alpha\in(0,\frac d2)\cap (0,1]$,
$\Lambda=\R^d$ is included (hence the cases studied in \cite{Kim}
are all covered by our results). Apart from this, the main novelty
of our applications is that (unlike in \cite{DRRW}) we can include
the case where $r_i<1$. Shortly speaking, the condition on $\Psi$ is
that the function
\[
  s\mapsto s\,\Psi(s)\;,\quad s\in\R,
\]
is (equal to or appropriately) comparable to a Young function $N$,
and $\Phi$ should be such that we can treat it as a kind of
perturbation (see conditions \textbf{(A1)}, \textbf{(A2)}, respectively,
in Section~\ref{sec:stoch-equations} below).
Then the appropriate choice of $V$ is
\[
  V := L_N\cap H\;,
\]
where $L_N$ is the Orlicz space determined by $N$ (cf.\ e.g.\
\cite{RZ}).
All details are contained in Section~\ref{sec:stoch-equations}.
We refer in particular to Theorem~\ref{thm-summary} there, which
summarizes the main results.

Finally, we mention that standardly by the variational approach one
also obtains information about the qualitative behaviour of solutions.
The corresponding results in our case are stated in Proposition~\ref{P2.2}
below.

\section{Monotone stochastic equations}\label{sec:mon-sdes}

We refer to \cite{KR} for the extensive literature on the subject
and in particular, to the pioneering work in the stochastic case
due to Pardoux \cite{P1, P2}. Let $G$ be a real separable Hilbert
space, $W_t$ a cylindrical Brownian motion on $G$, i.e. for an ONB
$\{g_1, g_2,\cdots\}$ of $G$, $W_t:=\big(B_t^i\, g_i\big)_{i\in
\N}$ for a sequence of independent one-dimensional Brownian
motions $\{B_t^i\}_{i\ge 1}$ on a complete filtered probability
spaces $(\OO, \F, \F_t; P)$. The filtration is assumed to be right
continuous. Next, let $(H, \<\, ,\,\>_H)$ be another real
separable Hilbert space. Let $V$ be a  reflexive Banach space such
that the embedding $V\subset H$ is dense and continuous, i.e. $V$
is dense in $H$ w.r.t. $\|\cdot\|_H$ and $\|v\|_H\le c\|v\|_V$ for
some constant $c>0$ and all $v\in V$. Let $V^*$ be the dual space
of $V$ and $_{V^*}\<\ ,\ \>_V$ denotes the corresponding
dualization. Identifying $H$ with its dual $H^*$ we have

$$V\subset H\equiv H^*\subset V^*$$ continuously and densely. Note
that, in particular, $V^*$ is then also separable, hence so is
$V$. Furthermore, $_{V^*}\<\ ,\ \>_V|_{H\times V}=\<\ ,\ \>_H.$
Let $\scr L_{HS}(G;H)$ denote the set of all real Hilbert-Schmidt
linear operators from $G$ to $H$, which is a real separable
Hilbert space under the inner product

$$\<B_1, B_2\>_{\L_{HS}}:= \sum_{i\ge 1} \<B_1 g_i, B_2 g_i\>_H.$$

We consider the following stochastic equation on $H$:

\beq\label{2.1} \d X_t = A(t,X_t)\d t + B(t,X_t)\d
W_t,\end{equation} where

$$A: [0,\infty)\times V \times  \OO\to V^*,\ \ \ B:
[0,\infty)\times V\times \OO\to \scr L_{HS}(G;H)$$ are
progressively measurable, i.e. for any $t\ge 0$, these mappings
restricted to $[0,t]\times V\times \OO$ are measurable w.r.t.
$\scr B([0,t])\times \scr B(V)\times \F_t$, where $\scr B(\cdot)$
is the Borel $\si$-field for a topological space. Below, writing
$A(t,v)$ we mean the mapping $\oo\mapsto A(t,v,\oo);$ analogously
for $B(t,v)$.

To solve  equation (\ref{2.1}), we need some assumptions. Let
$T>0$ be fixed. Let $(K,\|\cdot\|_K)$ be a real reflexive Banach
space such that $L^p([0,T]\times \OO\to V; \d t\times P)\subset
K\subset L^1([0,T]\times \OO\to V; \d t\times P)$ densely and
continuously for some $p\in (1,\infty)$. By Lemma \ref{L2} below,
its dual space is isometric to $(K^*,\|\cdot\|_{K^*})$, the
completion of $L^\infty([0,T]\times \OO\to V^*; \d t\times P)$
w.r.t.

\beq\label{*1}\|z^*\|_{K^*}:=\sup_{\|z\|_K\le 1}  \E\int_0^T \
_{V^*}\<z_t^*, z_t\>_V\,\d t.\end{equation} So, $K^*$ is reflexive
too and $K^*\subset L^{p/(p-1)}([0,T]\times\OO\to V^*;\d t\times
P)$ continuously and densely. Furthermore, $\ _{K^*}\<Y, z\>_K=
_{L^{p/(p-1)}}\<Y, z\>_{L^p}$ if $Y\in K^*$ and $z\in
L^p([0,T]\times \OO\to V; \d t\times P).$

Generalizing  the framework in \cite[Chapter 2]{KR}, we make the
following assumptions, where the first  is for the reference
spaces $K$ and $K^*$ and the remaining ones are for $A$ and $B$.

\begin{enumerate}
\item[{\bf(K)}]
  There
  exist a continuous function $R: V\to [0,\infty)$ with $R(x)=R(-x),
  x\in V,$ and $R(0)=0,$ and two locally bounded real functions
  $W_1, W_2$ on $[0,\infty)$ with $W_1(r),W_2(r)\to\infty$ as
  $r\to\infty$ such that:\newline
  \ \ (i) For any  sequence
  $z^{(n)}\in K, n\ge 1,\ \|z^{(n)}\|_K\to 0$ if and only if
  $\E\int_0^T R(z_t^{(n)})\d t\to 0$.\newline
  \ \ (ii) For any $z\in
  K$,\ $W_1\big(\E\int_0^T R(z_t)\d t\big)\le \|z\|_K\le W_2\big(
  \E\int_0^T R(z_t)\d t\big),$ where
  $W_1(\infty):=W_2(\infty):=\infty.$ \newline
  \ \  (iii) There exists a constant $C>0$ such that
  $R(x+y)\le C\bigl(R(2x)+R(2y))$ for all $x,y\in V.$\newline
  \ \ (iv) For
  any $h\in L^\infty([0,T]\times \OO;\d t\times P), z\in K$,
  we have $hz\in K$ and $\|hz\|_K\le \|h\|_\infty
  \|z\|_K$.
\item[{\bf(H1)}]
  Hermicontinuity: for any $u,v,x\in V, \oo\in \OO$ and any
  $t\in [0,T]$, the mapping
  $$\R\ni\ll\mapsto \ _{V^*}\<A(t,u+\ll v,\oo),x\>_V$$
  is continuous.
\item[{\bf (H2)}]
  Weak monotonicity: there exists a
  constant $c\in \R$ such that
  $$2\ _{V^*}\<A(\cdot,u)-A(\cdot,v), u-v\>_V + \|B(\cdot,u)-B(\cdot,v)\|_{\scr
    L_{HS}}^2\le c\|u-v\|_H^2\ \text{for\ all}\  u,v\in V$$ holds on
  $[0,T]\times \OO.$
\item[{\bf (H3)}]
  Coercivity: there exist
  constants $c_1,c_2>0$ and an $\F_t$-adapted process $f\in
  L^1([0,T]\times \OO;\d t\times P)$  such that $$2\, _{V^*}\<
  A(t,v), v\>_V+\|B(t, v)\|_{\L_{HS}}^2\le c_1\|v\|_H^2-c_2R(v)+
  f_t$$ holds on $\OO$ for all $v\in V$ and $t\in [0,T],$ where $R$
  is as in {\bf (K)}.
\item[{\bf (H4)}]
  There exists $c_3>0$ and an
  $\F_t$-adapted process $g\in L^1([0,T]\times \OO,\d t\times P)$
  such that
  $$|\, _{V^*}\<A(t, v),u\>_V|\le g_t+ c_3(R(v)+R(u))\ \text{on}\ \OO\
  \text{for\ all\ }t\in [0,T], u,v\in V,$$ where $R$ is as in {\bf (K)}.
\end{enumerate}

\paragraph{Remark 2.1.} (1) {\bf (H4)} together with {\bf
(K)}(ii) implies that for all $z\in K$ the map

$$\tt z\mapsto \E\int_0^T\, _{V^*}\<A(t,z_t),\tt
z_t\>_V\d t,\ \  \ \tt z\in K,$$ is in $K^*$ such that $K\ni
z\mapsto \|A(\cdot,z)\|_{K^*}$ is bounded on bounded sets in $K$.

(2) It easily follows from {\bf (H3)} and {\bf (H4)} that for all
$v\in V$ on $\OO$

$$ \|B(t,v)\|_{\L_{HS}}^2 \le
c_1\|v\|_H^2 + f_t+2g_t + 4 c_3 R(v),\ \ \ t\in [0,T].$$ In
particular, the function  $ K\in z\mapsto \E\int_0^T
\|B(t,z_t)\|_{\L_{HS}}^2\d t$ is by {\bf (K)}(ii) bounded on bounded
sets in $K$, which are also bounded in $L^2([0,T]\times \OO\to H; \d
t\times P)$.

(3) We observe that  $z\in K$ if and only if $z\in L^1([0,T]\times
\OO\to V; \d t\times P)$ with $\E\int_0^T
 R(z_t)\d t<\infty$. The  the necessity is trivial by {\bf (K)}(ii). On the other hand,
 for $z\in L^1([0,T]\times
\OO\to V; \d t\times P)$ with $\E\int_0^T
 R(z_t)\d t<\infty$ we have

 $$K \ni z^{(n)}:= z1_{\{|z|_V\le n\}}\to z\ \ \text{in\ }
 L^1([0,T]\times \OO\to V; \d t\times P).$$
 We claim that  $\{z^{(n)}\}$ is a Cauchy sequence in $K$ so that
 it also converges to $z$ in $K.$ Otherwise, there exist $\vv>0$ and a subsequence
 $n_k\to\infty$ such that

$$\sup_{m>n_k}\|z^{(n_k)}- z^{(m)}\|_K\ge 2\vv,\ \ \ k\ge 1.$$Then
for any $k\ge 1$

$$m_k:= \inf\{m>n_k: \|z^{(n_k)}- z^{(m)}\|_K\ge \vv\}<\infty.$$
Letting $\tt z^{(k)}= z^{(n_k)}-z^{(m_k)}$ we have $\|\tt
z^{(k)}\|_K\ge \vv$ for all $k\ge 1$, which is contradictive to
{\bf (K)}(i) since  $R(0)=0$ and the dominated convergence theorem
imply that

\[
\limsup_{k\to\infty}\E\int_0^T
 R(\tt z_t^{(k)})\d t\le \lim_{k\to\infty}\E\int_0^T R(z_t)1_{\{n_k<|z_t|_V\}}\d t
 = 0.
\]

(4) Let $Y\in K^*$, $z\in K$ and
$h\in L^\infty\bigl([0,T]\times\Omega;\mathrm dt\times P\bigr)$.
Let
$z^{(n)}\in L^p\bigl([0,T]\times\Omega\rightarrow V;\mathrm dt\times P\bigr)$,
$n\in\mathbb N$, such that $z^{(n)}\rightarrow z$ in $K$, hence
$h z^{(n)}\rightarrow hz$ in $K$ by \textbf{(K)}(iv).
Assume that $hY$ (which a priori is only in
$L^{p/(p-1)}\bigl([0,T]\times \Omega\to V^*;\mathrm dt\times P\bigr)$)
is in $L^\infty\bigl([0,T]\times \Omega\to V^*;\mathrm dt\times P\bigr)$.
Then
\begin{equation}\label{eq-2.2a}
  \pairing{K^*}{Y,hz}{K}
  = \E \int_0^T \pairing{V^*}{h_t Y_t,z_t}{V}\;\mathrm dt\;,
\end{equation}
since
\begin{align*}
  \pairing{K^*}{Y,hz}{K}
  &= \lim_{n\to\infty} \pairing{K^*}{Y,hz^{(n)}}{K}
  = \lim_{n\to\infty} \E \int_0^T \pairing{V^*}{Y_t,h_tz_t^{(n)}}{V}\;\mathrm dt\\
  &= \lim_{n\to\infty} \E \int_0^T \pairing{V^*}{h_t Y_t, z_t^{(n)}}{V}\;\mathrm dt
  = \E \int_0^T \pairing{V^*}{h_t Y_t,z_t}{V}\;\mathrm dt\;,
\end{align*}
since $z^{(n)}\to z$ in $K$, hence in
$L^1\bigl([0,T]\times \Omega\to V;\mathrm dt\times P\bigr)$.

(5) For $Y\in K^*$ we have
\begin{equation}\label{**}
  \ _{K^*}\<Y,z\>_K
  =\E\int_0^T\ _{V^*}\<Y_s, z_s\>_V\d s
  \le \E\int_0^T|\ _{V^*}\<Y_s, z_s\>_V|\d s
  \le \norm{Y}_{K^*}\, \|z\|_K,\ \ \ z\in K.
\end{equation}
Indeed, let $\xi:= _{V^*}\<Y,z\>_V$ and for $N\geq 1$
\begin{equation*}
  \begin{split}
    z^{(N)}
    &:= \text{sgn}(\xi) 1_{\{\|z\|_V\le N, \|Y\|_V*\le N\}}z
        \in L^p([0,T]\times \OO\to V; \d t\times P),\\
    Y^{(N)}
    &:= 1_{\{\|z\|_V\le N, \|Y\|_V*\le N\}}Y
        \in L^{p/(p-1)}([0,T]\times \OO\to V^*; \d t\times P).
  \end{split}
\end{equation*}
Then for all $N\geq 1$
\begin{equation*}
  \begin{split}
    &\E\int_0^T | \xi_t|1_{\{\|z\|_V\le N,
            \|Y\|_V*\le N\}}(t)\d t
    =\E\int_0^T \ _{V^*}\<Y_t^{(N)},
    z_t^{(N)}\>_V\d t =
    \ _{K^*}\<Y, z^{(N)}\>_K\\
    &\le \norm{Y}_{K^*}\, \|z^{(N)}\|_K
    \le \norm{Y}_{K^*}\,  \|z\|_K
  \end{split}
\end{equation*}
where we used \eqref{eq-2.2a} for the second equality
and \textbf{(K)}(iv) in the last step.
This implies \eqref{**} by letting $N\to \infty$.

(6)
By \cite[Proposition~26.4]{Zei90}, \textbf{(H1)} and
\textbf{(H2)} imply that for all $(t,\omega)\in[0,T]\times\Omega$
the map $u\mapsto A(t,u,\omega)$ is demicontinuous (i.e.\
$u_n\to u$ strongly in $V$ implies $A(t,u_n,\omega)\to A(t,u,\omega)$
weakly in $V^*$).
In particular, $A$ is continuous if $V$ is finite dimensional.
\\

We remark that the assumptions made in \cite[Chapter II]{KR} imply
the present ones by taking  $K:=L^p([0,T]\times\OO\to V; \d
t\times P)$ for some $p>1$, so that one has $K^*=
L^{p/(p-1)}([0,T]\times\OO\to V^*; \d t\times P),$ and
$R:=\|\cdot\|_V^p, W_1=W_2:= |\cdot|^{1/p}.$

\paragraph{Definition 2.1.} A   continuous adapted process
$\{X_t\}_{t\in [0,T]}$ on  $H$ is called a solution to
(\ref{2.1}), if $X\in L^2([0,T]\times\OO\to H, \d t\times P)$,
there exists a $\mathrm dt\times P$-version
$\bar X$ of an element in $K$ such that $X=\bar X\ \d t\times
P$-a.e., and $P$-a.s.

$$X_t= X_0 +\int_0^t A(s,\bar X_s)\d
s+\int_0^t B(s,\bar X_s)\d W_s,\  \text{for\ all\ } t\in [0,T].$$

\ \newline Since the Bochner integral on general Banach spaces is
only meaningful for pointwise (rather than a.e.) Banach
space-valued functions, in Definition 2.1 we have to choose a
$V$-valued progressively measurable version
of $X$ such that the right-hand side of the
above formula makes sense.

\beg{thm}\label{T2.1} Under ${\bf (K), (H1), (H2), (H3)}$ and {\bf
(H4)}, for any $X_0\in L^2(\OO\to H, \scr F_0; P),\ (\ref{2.1})$
has a unique solution and  the solution satisfies $\E\sup_{t\in
[0,T]}\|X_t\|_H^2<\infty$. Moreover; if $A(t,\cdot)(\oo)$ and
$B(t,\cdot)(\oo)$ are independent of  $t\in [0,T]$ and
$\oo\in\OO$, then the solution is a Markov process.
\end{thm}

The uniqueness of the solution can be easily proved using the It\^o
formula for $\|X_t\|_H^2$ presented in Theorem \ref{TA} in the
Appendix (see \cite[Theorem I.3.2]{KR} for a special case). In fact,
we have an even stronger statement formulated in the following
proposition. Furthermore, in the case where $A$ and $B$ are
independent of $t\in [0,T]$ and $\oo\in\OO$, the Markov property can
be easily proved by using the uniqueness as in \cite{KR}.

\beg{prp} \label{P2.2} Assume ${\bf (K), (H1), (H2), (H3)}$ and
{\bf (H4)}. Let $X$ and $Y$ be two solutions of $(\ref{2.1})$ with
$X_0,Y_0\in L^2(\OO\to H, \F_0; P)$. Let $c$ be the $($not
necessarily positive$)$ constant such that {\bf (H2)} holds. Then

\beq\label{R}\E\|X_t-Y_t\|_H^2\le \e^{ct} \E\|X_0-Y_0\|_H^2,\ \ \
t\in [0,T].\end{equation} Consequently,  if moreover $A$ and $B$
are independent of
 $t\in [0,T]$ and $\oo\in\OO$, the semigroup  $(P_t)_{t\in [0,T]}$ of the
corresponding Markov process is a Feller semigroup with

$$|P_t F(x)- P_t F(y)|\le \e^{ct/2}\text{Lip}(F)\|x-y\|_H,\ \ t\in
[0,T], x,y\in H,$$ where $F$ is an $H$-Lipschitz function with
Lip$(F)$ the Lipschitz constant. In particular, if  our
assumptions hold for each $T>0$ with $\ll:=-c>0$ independent of
$T$,  then $(P_t)$ has a unique invariant probability measure
$\mu$ with $\mu(\|\cdot\|_H^2)<\infty$ and $(P_t)$ converges
exponentially fast
 to $\mu$; more precisely,  for any $H$-Lipschitz function $F$,

$$|P_t F(x)-\mu(F)|^2\le \text{Lip}(F)^2\e^{-\ll t}
\mu(\|x-\cdot\|_H^2),\ \ \ x\in H.$$
\end{prp}

\beg{proof} By {\bf (H2)} and the It\^o formula (\ref{Ito}) in the
Appendix applied to $X_t-Y_t$,  we have
\begin{equation*}
  \begin{split}
    &\E\|X_t-Y_t\|_H^2-\E\|X_0-Y_0\|_H^2\\
    &= 2\E\int_0^t\big\{\ _{V^*}\<A(s,\bar X_s)-A(s,Y_s), \bar X_s-\bar Y_s\>_V
      +\|B(s,\bar X_s)-B(s,\bar Y_s)\|_{\L_{HS}}^2\big\}\d s\\
    &\le c\int_0^t\E\|X_s-Y_s\|_H^2\d s.
  \end{split}
\end{equation*}
This implies the first result immediately by Gronwall's lemma. The
remainder of the proof is the same as that of (3) and (4) in
\cite[Theorem 1.3]{DRRW}.\end{proof}

To prove the existence, we will construct a solution by the
classical Galerkin method of finite-dimensional approximations as
made in \cite{KR} (see also \cite{DRRW} for a special case).

Let $\{e_1,\cdots, e_n,\cdots\}\subset V$ be an ONB of $H$. Let
$H_n:=\text{span}\{e_1,\cdots, e_n\}, n\ge 1,$ and $P_n: V^*\to
H_n$ is defined by

\beq\label{Tr0}P_n y:= \sum_{i=1}^n \, _{V^*}\<y,e_i\>_V e_i,\ \ \
y\in V^*.\end{equation} Clearly, $P_n|_H$ is just the orthogonal
projection onto $H_n$ in $H$. Set

$$W_t^{(n)}:= \sum_{i=1}^n\<W_t,g_i\>_Gg_i=\sum_{i=1}^n B_t^i\,
g_i.$$ For each finite $n\ge 1$, we consider the following
stochastic equation on $H_n:$

\beq\label{F}\d\<X_t^{(n)}, e_j\>_H = \
_{V^*}\<A(t,X_t^{(n)}),e_j\>_V\d t+\<B(t,X_t^{(n)})\d W_t^{(n)},
\e_j\>_H,\ \ \ 1\le j\le n,\end{equation} where $X_0^{(n)}:= P_n
X_0.$ It is easily seen that we are in the situation covered by
\cite[Theorem 1.2]{K} which implies that (\ref{F}) has a unique
continuous strong solution.
 Let
\beq\label{R2}J:= L^2([0,T]\times \OO\to \scr L_{HS}(G;H); \d
t\times P).\end{equation}To construct the solution to (\ref{2.1}),
we need the following lemma.

\beg{lem} \label{L2.1} Under the assumptions in Theorem
$\ref{T2.1}$, we have
$$   \|X^{(n)}\|_K + \|A(\cdot, X^{(n)})\|_{K^*}+\sup_{t\in [0,T]}\E \|X_t^{(n)}\|_H^2
\le C$$ for some constant $C>0$ and all $n\ge 1.$
\end{lem}

\beg{proof} By It\^o's formula and {\bf (H3)}, we have

\beg{equation*}\beg{split} \d \|X_t^{(n)}\|_H^2 &= 2 \
_{V^*}\<A(t,X_t^{(n)}), X_t^{(n)}\>_V \d t +\|P_n B(t,
X_t^{(n)})\|_{\scr
L_{HS}(G;H)}^2\d t +\d M_t^{(n)}\\
&\le [-c_2 R(X_t^{(n)}) + c_1 \|X_t^{(n)}\|_H^2 + f_t]\d t+\d
M_t^{(n)}.\end{split}\end{equation*} for a local martingale
$M_t^{(n)}$. This implies

\beq\label{D*}\beg{split} &-\E\|X_0\|_H^2\le  \E
\e^{-c_1t}\|X_t^{(n)}\|_H^2 -\E\|X_0^{(n)}\|_H^2\\
&\le c_3 -c_4 \E\int_0^tR(X_s^{(n)})\d s,\ \ \ t\in
[0,T]\end{split}\end{equation} for some constants $c_3,c_4>0.$
Then the proof is completed by {\bf (K)} and {\bf
(H4)}.\end{proof}

\beg{lem} \label{L2} $K^*$ is isometric to the dual space $K'$ of
$K$.\end{lem}

\beg{proof} Since $L^\infty([0,T]\times \OO\to V^*;\d t\times P)$
is the dual space of $L^1([0,T]\times \OO\to V;\d t\times P)$, and
since the embedding $K\subset L^1([0,T]\times \OO\to V;\d t\times
P)$ is dense and continuous, the embedding

$$L^\infty([0,T]\times \OO\to V^*;\d
t\times P)\subset K'$$ is continuous too, but also dense by the
Hahn-Banach Theorem because $K$ is reflexive. Therefore, $K'=K^*$
since $K'$ is complete and contains $L^\infty([0,T]\times \OO\to
V^*;\d t\times P).$\end{proof}

\ \newline \emph{Proof of Theorem \ref{T2.1}.} We first recall
that the uniqueness of the solution is included in Proposition
\ref{P2.2}. Hence, when $A$ and $B$ are independent of the time
$t$ and $\oo\in\OO$, the Markov property follows immediately as in
the proof of \cite[Theorem II.2.4]{KR}. So, we only need to prove
the existence,  $\E\sup_{t\in [0,T]}\|X_t\|_H^2<\infty$  and the
continuity of $X_t$ in $H$.

(a)  By the reflexivity of $K$ and Lemmas \ref{L2.1}, \ref{L2},
and Remark 2.1, we have, for a subsequence $n_k\to\infty$:
\begin{enumerate}
    \item [(i)] $X^{(n_k)}\to \bar X \ \text{weakly\ in}\ K \ \text{and \ weakly\
in}\ L^2([0,T]\times \OO\to H; \d t\times P);$
 \item[(ii)] $A(\cdot, X^{(n_k)})\to Y$\ weakly in $K^*$ (hence weakly in $L^{p/(p-1)}
 ([0,T]\times\OO\to V^*;\d t\times P))$;\item[(iii)]  $B(\cdot, X^{(n_k)})\to
 Z$ weakly in $J$ and hence $\int_0^\cdot B(s,X_s^{(n_k)})\d
 W_s^{(n_k)} \to \int_0^\cdot Z_s\d W_s$ weakly in $L^\infty([0,T]\to L^2(\OO\to H; P); \d t)\
 ($equipped with the supremum norm).
\end{enumerate} Here the second part in (iii) follows since also $B(\cdot, X^{(n_k)})
\tt P_{n_k}\to Z$ weakly in $J$,  where $\tt P_n$ is the
orthogonal projection onto span$\{g_1,\cdots, g_n\}$ in $G$, since

$$\int_0^\cdot B(s, X_s ^{(n_k)})\d W_s^{(n_k)} =  \int_0^\cdot B(s, X_s^{(n_k)})\tt P_{n_k}\d W_s$$
and since a bounded linear operator between two Banach spaces is
trivially weakly continuous. Since the approximants are
progressively measurable, so are $\bar X, Y$ and $Z$.

Thus, by the definition of $X^{(n)}$ we have (since $V$ is
separable) that

\beq\label{2.2} \ _{V^*}\<\bar X_t,e\>_V =\<X_0, e\>_H +\int_0^t \
_{V^*}\<Y_s, e\>_V\d s +\int_0^t \<Z_s\d W_s, e\>_H,\ \text{for\ all
\ } e\in V\  \text{a.e.-}\d t\times P.\end{equation} Defining

\beq\label{2.3} X_t:= X_0 +\int_0^t Y_s\d s+\int_0^t Z_s\d W_s,\ \
\ t\in [0,T],\end{equation} we have $X=\bar X$ $\d t\times P$-a.e.
 We note here, that since $Y\in K^*\subset L^{p/(p-1)}([0,T]\times \OO\to V^*;
 \d t\times P)$, its integral in (\ref{2.3}) always exists as
 a $V^*$-valued Bochner integral and is continuous in $t$.
Therefore,  $X$ is a $V^*$-valued continuous adapted process.

Hence Theorem \ref{TA} applies to $X$ in
(\ref{2.3}), so $X$ is continuous in $H$ and $\E \sup_{t\le T}
\|X_t\|_H^2<\infty.$

Thus, it remains to verify that

\beq\label{2.4} B(\cdot, \bar X)= Z,\ \ \ A(\cdot, \bar X)= Y,\ \
\ \text{a.e.-}\d t\times P.\end{equation} To this end, we first
note that for any nonnegative $\psi\in L^\infty([0,T],\d t)$ it
follows from (i)  that

\beg{equation*}\beg{split}\E\int_0^T \psi(t) \|\bar X_t\|_H^2\d
t&=\lim_{k\to\infty} \E\int_0^T\<\psi(t)\bar X_t, X_t^{(n_k)}\>_H\d t\\
&\le \bigg( \E\int_0^T \psi(t) \|\bar X_t\|_H^2\d
t\bigg)^{1/2}\liminf_{k\to\infty} \bigg( \E\int_0^T \psi(t) \|
X_t^{(n_k)}\|_H^2\d
t\bigg)^{1/2}<\infty.\end{split}\end{equation*} Since $X=\bar X\
\d t\times P$-a.e., this implies
 \beq\label{2.5} \liminf_{k\to\infty} \E\int_0^T\psi(t)
\|X_t^{(n_k)}\|_H^2\d t \ge \E \int_0^T \psi(t) \|X_t\|_H^2\d
t.\end{equation} By (\ref{2.3}) and the It\^o formula (\ref{Ito})
in the Appendix, using the elementary product rule we obtain that

\beq\label{2.6} \E \e^{-ct} \|X_t\|_H^2 - \E\|X_0\|_H^2
=\E\int_0^t\e^{-cs}\big\{2\ _{V^*}\<Y_s, \bar X_s\>_V
+\|Z_s\|_{\L_{HS}}^2-c\|X_s\|_H^2\big\}\d
s.\end{equation}Furthermore, for any $\phi\in K\cap
L^2([0,T]\times \OO\to H; \d t\times P)$,

\beq\label{2.7}\beg{split}& \E \e^{-ct}\|X_t^{(n_k)}\|_H^2
-\E\|X_0^{(n_k)}\|_H^2\\
&=\E \int_0^t\e^{-cs}\big\{2\ _{V^*}\<A(s,X_s^{(n_k)}),
X_s^{(n_k)}\>_V +\|
P_{n_k}B(s,X_s^{(n_k)})\tt P_{n_k}\|_{\L_{HS}}^2-c\|X_s^{(n_k)}\|_H^2\big\}\d s\\
&\le \E \int_0^t\e^{-cs}\big\{2\ _{V^*}\<A(s,X_s^{(n_k)}),
X_s^{(n_k)}\>_V +\|
B(s,X_s^{(n_k)})\|_{\L_{HS}}^2-c\|X_s^{(n_k)}\|_H^2\big\}\d
s.\\
&=\E\int_0^t \e^{-cs}\Big\{2 \
_{V^*}\<A(s,X_s^{(n_k)})-A(s,\phi_s), X_s^{(n_k)}-\phi_s\>_V\\
&\quad+\|B(s,X_s^{(n_k)})-B(s,\phi_s)\|_{\L_{HS}}^2-c\|X_s^{(n_k)}-\phi_s\|_H^2\Big\}\d
s\\
&\quad + \E\int_0^t\e^{-cs}\Big\{2\ _{V^*}\<A(s,\phi_s),
X_s^{(n_k)}\>_V + 2\ _{V^*}\<A(s,X_s^{(n_k)})-A(s,\phi_s),
\phi_s\>_V\\
&\quad -\| B(s,\phi_s)\|_{\L_{HS}}^2 + 2
\<B(s,X^{(n_k)}_s),B(s,\phi_s)\>_{\L_{HS}} -2c\<X_s^{(n_k)},
\phi_s\>_H+c\|\phi_s\|_H^2\Big\}\d s.
\end{split}\end{equation}
Note that by {\bf (H2)} the first of the two summands above is
negative. Hence by letting $k\to\infty$ we conclude by (i)--(iii),
Fubini's theorem, Remark 2.1,  (\ref{**}) and (\ref{2.5}) that for
every nonnegative $\psi\in L^\infty([0,T];\d t)$

\beg{equation*}\beg{split}&\E
\int_0^T\psi(t)(\e^{-ct}\|X_t\|_H^2-\|X_0\|_H^2)\d t \\
&\le \E\int_0^T\psi(t)\bigg(\int_0^t\e^{-cs}\Big\{2\
_{V^*}\<A(s,\phi_s), \bar X_s\>_V + 2\ _{V^*}\<Y_s-A(s,\phi_s),
\phi_s\>_V\\
&\quad -\| B(s,\phi_s)\|_{\L_{HS}}^2 + 2
\<Z_s,B(s,\phi_s)\>_{\L_{HS}} -2c\<X_s,
\phi_s\>_H+c\|\phi_s\|_H^2\Big\}\d s\bigg)\d t.
\end{split}\end{equation*}
Inserting (\ref{2.6}) for the left hand side and rearranging as
above we arrive at

\beq\label{2.8} \beg{split} 0\ge &\E\int_0^T\psi(t) \bigg(
\int_0^t\e^{-cs}\big\{
2 \ _{V^*}\<Y_s-A(s,\phi_s), \bar X_s-\phi_s\>_V\\
& + \|B(s,\phi_s)-Z_s\|_{\L_{HS}}^2-c\|X_s-\phi_s\|_H^2\big\}\d
s\bigg)\d t.
\end{split}\end{equation}Taking $\phi=\bar X$ we obtain from (\ref{2.8}) that $Z= B(\cdot, \bar X)$.
Finally, first applying (\ref{2.8}) to $\phi= \bar X-\vv \tt\phi\,
e$ for $\vv>0$ and $\tt \phi\in L^\infty([0,T]\times \OO; \d
t\times P), e\in V$, then dividing both sides by $\vv$ and letting
$\vv\to 0$, by the dominated convergence theorem, the
hemicontinuity of $A$, {\bf (K)} and {\bf (H4)}, we obtain

$$0\ge \E\int_0^T \psi(t) \bigg(\int_0^t \e^{-cs}\tt \phi_s\
_{V^*}\<Y_s-A(s,\bar X_s),e\>_V\d s\bigg)\d t.$$ By the
arbitrariness of $\psi$ and $\tt\phi$, we conclude that $Y=
A(\cdot, \bar X).$ This completes the proof.\qed

\section{Stochastic generalized porous medium and fast diffusion
equations}\label{sec:stoch-equations}

In this section we shall discuss applications.
Let $(E,\scr B, \m)$ be a $\si$-finite measure space
with countably generated $\sigma$-algebra $\scr B$.
Let
$(L,\D(L))$ be a negative definite self-adjoint operator on
$L^2(\m)$ with Ker$(L)=\{0\}$. We shall use $\<\ ,\ \>$ and
$\|\cdot\|$ to stand for the inner product and the norm in
$L^2(\m)$ respectively, we also denote
$\<f,g\>:=\m(fg) := \int fg \d \m$ for any
two functions $f,g$ such that $fg\in L^1(\m)$.
Consider the quadratic form $(\scr E, D(\scr E))$ on $L^2(\m )$
associated with $(L,\scr D(L))$, i.e.
\[
  D(\scr E) := \scr D\bigl(\sqrt{-L}\bigr)
\]
and
\[
  \scr E (u,v)
  := \m \bigl(\sqrt{-L} u \, \sqrt{-L} v\bigr)\;;
  \quad u,v\in D(\scr E).
\]
Let $\scr F_e$ be the abstract completion of $D(\scr E)$ with respect
to the norm
\[
  \|\,\cdot\,\|_{\scr F_e} := \scr E(\, \cdot \, ,\, \cdot \, )^{\frac 12}\;,
\]
and let $\scr F_e^*$ be its dual space.
Note that both $(\scr F_e,\scr E)$ and $\scr F_e^*$, with the
inner product induced by the Riesz isomorphism, are
Hilbert spaces.
Now we define
\begin{equation}\label{eq-3.0}
  H := \scr F_e^* \;,
  \quad \< \, \cdot \, , \, \cdot \, \>_H
        := \< \, \cdot \, , \, \cdot \,  \>_{\scr F_e^*}\;,
\end{equation}
i.e.\ $\scr F_e^*$ will be the state space of our SDE
\eqref{2.1}.
In order to make our general results from the previous
section work, we shall use the further assumption, that
$(\scr E, D(\scr E))$ is a transient Dirichlet space in
the sense of
\cite[Section~1.5]{FOT94}.
$(\scr E,\scr F_e)$, where $\scr E$ also denotes the extension
of $\scr E$ to $\scr F_e$, is called \emph{extended Dirichlet
space} in \cite{FOT94}, from which we also adopt the notation.
In this case, $\scr F_e^*$ is also called the corresponding
\emph{Green space}.
Since we do not assume the reader to be familiar with all
these notions, we shall first formulate some abstract
conditions on $(\scr E, D(\scr E))$ (i.e.\ on $(L,D(L))$)
such that our proofs below work.
Then, partly on the basis of results from \cite{FOT94},
we shall prove that these abstract assumptions hold
if $(\scr E, D(\scr E))$ is a transient Dirichlet space.
Furthermore, we shall briefly describe whole classes
of concrete examples.

From now on, we are going to assume that the following condition
holds:
\begin{enumerate}
\item[\textbf{(L1)}]
  There exists a strictly positive $g\in L^1(\m)\cap L^\infty(\m)$
  such that $\scr F_e\subset L^1(g\cdot \m)$ continuously.
\end{enumerate}
Our space $V$ will be heuristically given as $V := H\cap L_N$, where
$N$ is a nice Young function and $L_N$ the corresponding Orlicz
space. More precisely, let $N\in C(\R)$ be a Young function, i.e. a
nonnegative, continuous, convex and even function such that $N(s)=0$
if and only if $s=0$, and

$$\lim_{s\to 0} \ff{N(s)}s=0,\ \ \ \lim_{s\to\infty}\ff{N(s)}s
=\infty.$$ For any function $f$ on $E$ with $\m(N(\aa f))<\infty$
for some $\aa>0$, define

$$\|f\|_{L_N}:=\inf\{\ll\ge 0: \m(N(f/\ll))\le 1\}.$$
Then the space $(L_N , \| \, \cdot \, \|_{L_N})$, where

$$L_N(\m) := \{f: \|f\|_N<\infty\} , $$
is a  real separable Banach space, which is called the Orlicz
space induced by the Young-function $N$ (cf. \cite[Proposition
1.2.4]{RZ}). There is an equivalent norm defined by using the dual
function:

$$N^*(s):= \sup\{r|s|- N(r): r\ge 0\},\ \ \ s\in\R,$$
which is once again a Young function. More precisely, letting

$$\|f\|_{(N)}:=\sup\{\<f,g\>: \m(N^*(g))\le 1\},$$
one  has (see \cite[Theorem 1.2.8 (ii)]{RZ})

\beq\label{*0}
  \|\cdot\|_{L_N}\le \|\cdot \|_{(N)}\le 2
  \|\cdot\|_{L_N}.
\end{equation}
The function $N$ is called
$\DD_2$-regular, if there exists constant $c>0$ such that

$$N(2s)\le c\big(N(s) + 1_{\{\m(E)<\infty\}}\big),\ \ \ s\in\R.$$
In this paper we assume that $N$ and $N^*$ are $\DD_2$-regular. By
\cite[Proposition 1.2.11(iii) and Theorem 1.2.13]{RZ}, $L_N(\m)$
and $L_{N^*}(\m)$ are dual spaces of each other
with dualization given by $\<f,g\>=\m (fg)$, $f\in L_N$, $g\in L_{N^*}$.
Hence they are
reflexive. By the $\DD_2$-regularity, $f\in L_N(\m)$ if and only
if $\m(N(f))<\infty$.

Now let us give a precise definition of $V$. Define $V:= H\cap
L_N(\m)$ in the following sense:
\begin{equation}
  V
  := \bigl\{
       u\in L_N(\m )
     \bigm|
       \exists\, c\in(0,\infty) \text{ such that }
       \m (uv)\leq c\,\|v\|_{\scr F_e} \; \forall\, v\in\scr F_e\cap L_{N^*}
     \bigr\}
\end{equation}
equipped with the norm
\[
  \|u\|_V
  := \|u\|_{L_N} + \|u\|_{\scr F_e^*}
  = \|u\|_{L_N} + \|u\|_H\;.
\]
In order that $V$ becomes a subset of $H$ and
$(V,\|\, \cdot \, \|_V)$ a Banach space, we need to assume that
$\scr F_e \cap L_{N^*}$ is a dense subset of $\scr F_e$.
For later use, we even make the following stronger assumption:
\begin{enumerate}
\item[\textbf{(N1)}]
  $\scr F_e\cap L_{N^*}$ is a dense subset of both $\scr F_e$
  and $L_{N^*}$.
\end{enumerate}
By \textbf{(N1)}, $V$ can be considered as a subset of $H=\scr F_e^*$
by identifying $u\in V$ with the map
$\bar u : \scr F_e\cap L_{N^*}\to\mathbb R$ defined by
\[
  \bar u(v)
  := \m (uv)\;,
  \quad v\in \scr F_e \cap L_{N^*}\;.
\]
Then obviously, $V\subset H$ continuously, and it is easy
to see that $V$ is complete with respect to $\| \, \cdot \, \|_V$.
The density of $V$ in $H$ is, however, not clear.
Therefore, we assume:
\begin{enumerate}
\item[\textbf{(N2)}]
  $V$ is a dense subset of both $H=\scr F_e^*$ and $L_N$.
\end{enumerate}
The second part of \textbf{(N2)} we shall only need later.
As mentioned above, we shall later prove that \textbf{(L1)},
\textbf{(N1)} and \textbf{(N2)} always hold if
$(\scr E, D(\scr E))$ is a transient Dirichlet space.
Let $V^*$ be the dual space of $V$, so using $H\equiv H^*$
we have
\[
  V \subset H \subset V^*
  \quad \text{continuously and densely.}
\]
Note that, since $V$ is complete, the map
\[
  V\ni u \,\mapsto\, (u,\m(u\,\cdot\,)) \in L_N\times\scr F_e^*
\]
is an isomorphism from
$V$ to a closed subspace of $L_N\times \scr F_e^*$
which is reflexive.
So, $V$ itself is reflexive.

\begin{prp}\label{prp-3.1}
  Assume that $(\scr E,D(\scr E))$ is a transient Dirichlet space.
  Then conditions \textbf{\upshape (L1)}, \textbf{\upshape (N1)} and
  \textbf{\upshape (N2)} hold.
\end{prp}

With respect to the length of this paper we do not recall all
necessary definitions and notions here, but refer to
\cite[Section~1.5]{FOT94}.
We shall stick exactly to the terminology and notation introduced there.

Before we prove Proposition~\ref{prp-3.1}, we need the following
fact on $\Delta_2$-regular Young functions.

\begin{lem}\label{lem-3.2}
  Let $N$ be $\Delta_2$-regular.
  Then:
  \begin{enumerate}
  \item[(i)]
    There exists $q\in (2,\infty)$ such that
    \begin{equation}\label{eq:delta-reg}
      N(rs)
      \leq r^q \bigl( N(s) + 2\cdot 1_{\{\m(E)<\infty\}} \bigr)\;,
      \quad r\geq 2, \; s\geq 0.
    \end{equation}
  \item[(ii)]
    $L^1(\m)\cap L^q(\m)\subset L_N(\m)$ continuously,
    where the intersection is equipped with the norm
    $\| \; \|_1 + \| \; \|_q$ and $q$ is as in (i).
  \end{enumerate}
\end{lem}

\begin{proof} (i)
    Let $n\geq 1$ be such that $r\in[2^n,2^{n+1})$. If $N$ is
    $\Delta_2$-regular with constant $C>2$, then for
    $p_1 := \log C / \log 2$
    \begin{align*}
      N(rs)
      &\leq N(2^{n+1}s)
      \leq C^{n+1} \bigl( N(s) + \sum_{i=0}^n C^{-i}\cdot 1_{\{\m(E)<\infty\}} \bigr)
      \le 2^{p_1(n+1)} \bigl( N(s) + 2\cdot 1_{\{\m(E)<\infty\}} \bigr)\\
      &\leq r^{p_1(n+1)/n} \bigl( N(s) + 2\cdot 1_{\{\m(E)<\infty\}} \bigr)\;.
    \end{align*}
    Thus, \eqref{eq:delta-reg} holds by taking $q := 2p_1>2$.

  (ii)    We have for all $\lambda>0$ and $f\in L^1(\m)\cap L^q(\m)$
    \begin{align*}
      \m\left(N\biggl(\frac{|f|}{\lambda}\biggr)\right)
     &= \m\left( 1_{\{\frac{|f|}{\lambda}\geq 2\}}
                 N\biggl(\frac{|f|}{\lambda}\biggr)\right)
        + \m\left( 1_{\{\frac{|f|}{\lambda}< 2\}}
                   N\biggl(\frac{|f|}{\lambda}\biggr)\right)\\
     &\leq \frac{N(1)+2}{\lambda^q}\,\m(|f|^q)
           + \sup_{0\leq s\leq 2}\frac{N(s)}{s}\cdot \frac 1\lambda\,\m(|f|)
    \end{align*}
    where we used (i) in the last step.
    For
    \[
      \lambda
      := \bigl(2\m(|f|^q)(N(1)+2)\bigr)^{\frac 1q}
         + 2\m(|f|)\sup_{0\leq s\leq 2}\frac{N(s)}{s}\;,
    \]
    the right hand side is less than $1$.
    Hence, the assertion follows.
\end{proof}

\begin{proof}[Proof of Proposition~\ref{prp-3.1}]
  We first note that \textbf{(L1)} holds by
  \cite[Theorem~1.5.3 ($\alpha$) and ($\beta$)]{FOT94}
  (see also \cite[Lemma~1.5.5]{FOT94}).

  Now let us prove \textbf{(N1)}.
  We first note that $L^1(\m)\cap L^\infty(\m)$ is a dense subset
  of $L_{N^*}(\m)$.
  This follows from Lemma~\ref{lem-3.2}(ii) applied to $L_{N^*}$
  and since for any $f\in L_N$ such that
  \[
    \m(vf)=0
    \quad\text{for all $v\in L^1(\m)\cap L^\infty(\m)$}
  \]
  it follows that $f=0$.
  So, to show that $\scr F_e\cap L_{N^*}$ is dense in $L_{N^*}$,
  by Lemma~\ref{lem-3.2}(ii) it suffices to show that
  \begin{equation}\label{eq-3.4}
    D(\scr E) \cap L^1(\m)\cap L^\infty(\m)
    \quad\text{is dense in $L^1(\m)\cap L^q(\m)$,}
  \end{equation}
  where $q$ is as in Lemma~\ref{lem-3.2}(ii) with $N^*$
  replacing $N$.
  But \eqref{eq-3.4} is a well-known fact about Dirichlet spaces.
  To show that $\scr F_e \cap L_{N^*}$ is dense in $\scr F_e$
  it suffices to show that
  \begin{equation}\label{eq-3.5}
    L^1(\m)\cap L^\infty(\m)\cap D(\scr E)
    \quad\text{is dense in $D(\scr E)$}
  \end{equation}
  with respect to the norm
  $\scr E_1(\,\cdot\,,\,\cdot\,)^{\frac 12}
   := \bigl(\scr E(\,\cdot\,,\,\cdot\,)+\|\,\cdot\,\|_2\bigr)^{1/2}$.
  So, let $u\in D(\scr E)$ such that
  \[
    \scr E_1(v,u)
    = 0
    \quad\text{for all $v\in L^1(\m)\cap L^\infty(\m)\cap D(\scr E)$.}
  \]
  Then for $f\in L^1(\m)\cap L^2(\m)$, since
  $G_1 f\in L^1(\m) \cap L^\infty(\m)\cap D(\scr E)$, where
  $(G_\lambda)_{\lambda>0}$ is the resolvent associated to
  $(\scr E, D(\scr E))$,
  \[
    \m (fu)
    = \scr E_1 (G_1 f, u)
    = 0\;.
  \]
  Hence $u=0$, because $L^1(\m)\cap L^\infty(\m)$
  is dense in $L^2(\m)$.

  Now let us show \textbf{(N2)}.
  Let $g$ be as in \textbf{(L1)} and consider the set
  \[
    \scr G
    := \bigl\{ h\cdot g\bigm| h\in L^\infty(\m)\bigr\}\;.
  \]
  Then $\scr G\subset L^1(\m)\cap L^\infty(\m)\subset L_N(\m)$.
  Since $L^1(\m)\cap L^\infty(\m)$ is dense in $L_N(\m)$,
  it follows that $\scr G$ is dense in $L_N(\m)$.
  Furthermore, by
  \cite[Theorem~1.5.4]{FOT94} for every $f\in \scr G$
  there exists $Gf\in\scr F_e$ such that
  \[
    \m(f v)
    = \scr E(Gf,v)
    \quad\text{for all $v\in\scr F_e$.}
  \]
  Hence $\scr G\subset V$, so $V$ is dense in $L_N(\m)$.

  Now we show that $\scr G$ is dense in $H=\scr F_e^*$.
  So, let $v\in\scr F_e$ such that
  \[
    \m (fv)=0
    \quad\text{for all $f\in \scr G$.}
  \]
  Then
  \[
    \m(hvg) = 0
    \quad\text{for all $h\in L^\infty(\m)$,}
  \]
  and hence, since $vg\in L^1(\m)$, $vg=0$ $\m$-a.e.
  Since $g$ is strictly positive, it follows that $v=0$
  $\m$-a.e.\ and \textbf{(N2)} is proved.
\end{proof}

\begin{exa}\label{exa-3.3}
  There are plenty of examples of transient Dirichlet spaces
  described in the literature with $E$ being e.g.\ a
  manifold or a fractal.
  Here, we shall briefly refer to cases where $E:= D\subset \mathbb R^d$
  and $\m$ is the Lebesgue measure and which are presented in detail
  in
  \cite{FOT94} (see Examples 1.5.1--1.5.3).
  For example, if $L$ is the Friedrichs extension of a symmetric
  uniformly elliptic operator of second order with Dirichlet
  boundary conditions on an open domain $D\subset \mathbb R^d$,
  then all the above applies.
  But there are also examples with Neumann boundary conditions,
  as e.g.\ the Laplacian on $D:=$ half space in $\mathbb R^d$.
  Furthermore, for $D:=\mathbb R^d$ we can take $L=(-\Delta)^\alpha$ with its
  standard domain, if $\alpha\in(0,\frac d2)\cap (0,1]$.
  For details we refer to \cite{FOT94}.
\end{exa}

Now let us return to our general situation described at
the beginning of this section, i.e.\ the only conditions on
$L$ and the $\Delta_2$-regular (dual) Young functions $N$
and $N^*$ are \textbf{(L1)}, \textbf{(N1)} and \textbf{(N2)}.
The following is then standard and more a question of
notation than contents.
Nevertheless, we include a short proof.

\begin{lem}\label{lem-3.3}
  \begin{enumerate}
  \item[(i)]
    The map $\bar L:\scr F_e \to \scr F_e^*$ defined by
    \begin{equation}\label{eq-3.6}
      \bar Lv := -\scr E (v,\,\cdot\,)\;,\quad v\in \scr F_e
    \end{equation}
    (i.e.\ the Riesz isomorphism of $\scr F_e$ and $\scr F_e^*$
    multiplied by $(-1)$)
    is the unique continuous linear extension of the map
    \begin{equation}\label{eq-3.7}
      D(L)\ni v \;\mapsto\; \m(Lv\;\cdot\,)\in\scr F_e^*\;.
    \end{equation}
  \item[(ii)]
    Let $v\in \scr F_e\cap L_{N^*}$, $u\in\scr F_e^*\cap L_N=V$.
    Then
    \[
      \< \bar Lv, \bar u \>_{\scr F_e^*}
      = -\m(vu)\;.
    \]
  \item[(iii)]
    The map $\bar L : \scr F_e\cap L_{N^*}\to\scr F_e^*\subset V^*$
    has a unique continuous linear extension (again denoted by)
    $\bar L : L_{N^*}\to V^*$
    and this extension satisfies
    \begin{equation}\label{eq-3.8}
      \sideset{_{V^*}}{_{V}}{\mathop{\<\bar Lv, u\>}}
      = -\m(vu)
      \quad\text{for all $v\in L_{N^*}$, $u\in V$.}
    \end{equation}
  \end{enumerate}
\end{lem}

\begin{proof} (i)     For all $v\in D(L)$, $w\in D(\scr E)$ we have
    \[
      \m(Lv\,w) = -\scr E(v,w)\;.
    \]
    By the density of $D(\scr E)$ in its completion $\scr F_e$
    it follows that the linear map in \eqref{eq-3.7} really takes
    values in $\scr F_e^*$ and that it is continuous as a map
    from $\scr F_e$ with domain $D(L)$ with values in $\scr F_e^*$.
    Since $D(L)$ is dense in $(D(\scr E), \scr E_1)$, hence in
    $\scr F_e$, the assertion follows.

  (ii)    By Riesz's representation theorem there exists $Gu\in\scr F_e$
    such that $\bar u \;(:= \m(u\;\cdot\,))\;=\scr E(Gu,\,\cdot\,)$.
    Hence,
    \begin{align*}
      \<\bar Lv, \bar u\>_{\scr F_e^*}
      &= -\bigl\< \scr E(v,\,\cdot\,), \scr E(Gu,\,\cdot\,)\bigr\>_{\scr F_e^*}
       = -\scr E(v,Gu)
       = - \m(vu)\;.
    \end{align*}

  (iii)    Let $v\in\scr F_e\cap L_{N^*}$, $u\in\scr F_e^*\cap L_N = V$.
    Then by (ii)
    \begin{equation}\label{eq-3.9}
      \sideset{_{V^*}}{_{V}}{\mathop{\<\bar Lv, u\>}}
      = \< \bar Lv,\bar u \>_{\scr F_e^*}
      = -\m(vu)\;,
    \end{equation}
    hence,
    \[
      \bigl| \sideset{_{V^*}}{_{V}}{\mathop{\<\bar Lv, u\>}} \bigr|
      \leq c\,\|v\|_{L_{N^*}}\,\|u\|_{L_N}
      \leq c\, \|v\|_{L_{N^*}}\, \|u\|_V
    \]
    for some positive constant $c$ independent of $u$ and $v$.
    By \textbf{(N2)} it follows that
    \[
      \|\bar L v\|_{V^*}
      \leq c\, \|v\|_{L_{N^*}}
      \quad\text{for all $v\in\scr F_e\cap L_{N^*}$.}
    \]
    Hence, by the second part ot \textbf{(N1)}, the desired
    extension $\bar L : L_{N^*}\to V^*$ exists and
    \eqref{eq-3.9} extends to \eqref{eq-3.8}.
\end{proof}

For simplicity, we write $L$ instead of $\bar L$ and $u$
instead of $\bar u\in V$ below, hence, consider $V$ as
a subset of $H=\scr F_e^*$, hence of $V^*$ in particular.

To define the nonlinear operators $A$ and $B$ in \eqref{2.1},
let
\begin{equation}\label{eq-3.9a}
  \Psi,\Phi: [0,\infty)\times \R\times \OO\to \R
\end{equation}
be progressively measurable, i.e. for any $t\ge 0$, restricted to
$[0,t]\times \R\times \OO$ they are measurable w.r.t. $\scr
B([0,t])\times\scr B(\R)\times \F_t.$ We assume that for fixed
$(t,\oo)\in [0,\infty)\times \OO,$ $\Psi(t,\cdot)(\oo)$ and
$\Phi(t,\cdot)(\oo)$ are continuous.

Finally, let
\begin{gather}\label{eq-3.9b}
\begin{tabular}{p{13cm}}
$B: [0,T]\times V\times \OO\to \L_{HS}(G;H)$ be progressively
measurable satisfying \textbf{(H2)} with $A\equiv 0$ (i.e., $B$ is
Lipschitz with respect to $H$-norm, uniformly in
$(t,\omega)\in[0,T]\times \Omega$).\end{tabular}
\end{gather}
We distinguish two sets of conditions on $\Psi$ and $\Phi$:
\begin{enumerate}
\item[\textbf{(A1)}]
  $\Phi(t,s)=h_ts$, $t\in[0,T]$, $s\in \mathbb R$, for some
  $\scr F_t$-adapted
  $h\in L^\infty \bigl([0,T]\times\Omega;\mathrm dt\times P\bigr)$
  and there exist $\Delta_2$-regular dual Young
  functions $N$ and $N^*$, a nonnegative $\scr F_t$-adapted process
  $f\in L^1\bigl([0,T]\times\Omega;\d t\times P\bigr)$
  and a constant $c\geq 1$, such that for all
  $s,s_1,s_2\in \mathbb R$ on $[0,T]\times \Omega$
  \begin{gather}
  \tag{$\Psi 1$}
    (s_2-s_1)\bigl(\Psi(\,\cdot\,,s_2)-\Psi(\,\cdot\,, s_1)\bigr)
    \geq 0\;. \\
  \tag{$\Psi 2$}
    s\Psi(\,\cdot\,,s)
    \geq N(s)- 1_{\{\m(E)<\infty\}} \cdot f \;. \\
  \tag{$\Psi 3$}
    s\Psi(\,\cdot\,,s)
    \leq c\bigl( N(s)+ 1_{\{\m(E)<\infty\}} \cdot f \bigr) \;.\\
  \tag{$\Psi 4$}
    N^*\bigl(\Psi(\argdot,0)\bigr)\, 1_{\{\m(E)<\infty\}}
    \in L^1\bigl( [0,T]\times\Omega;\mathrm dt\times P \bigr)\;.
  \end{gather}
\end{enumerate}
If $\Phi$ is not just the identity times $h$ as above,
we need to restrict to the Young function considered in
Example~\ref{exa-3.5}(i), i.e.
\begin{equation}\label{eq-3.10}
  N(s)
  := \sum_{i=1}^m \varepsilon_i |s|^{r_i+1}\;,
  \quad s\in \mathbb R,
\end{equation}
for some pairwise distinct $r_1,\dots,r_m>0$ and
$\varepsilon_1,\dots,\varepsilon_m>0$, and consider the
following condition:
\begin{enumerate}
\item[\textbf{(A2)}]
  There exist $r_1,\dots,r_m>0$ such that
  for $N$ as in \eqref{eq-3.10} the set
  $L\bigl( D(L)\cap L_N(\m) \bigr)$
  is dense in $L^{r_i+1}(\m)$ and
  $L^{-1}:L^{r_i+1}(\m)\to L^{r_i+1}(\m)$, $1\leq i\leq m$,
  is bounded.
  Furthermore, there exist a nonnegative $\scr F_t$-adapted process
  $f\in L^1\bigl([0,T]\times\Omega, \mathrm dt\times P\bigr)$
  and  constants $c\geq 0$, $\delta_1,\dots,\delta_m>0$
  such that for all $s,s_1,s_2\in\mathbb R$
  we have on $[0,T]\times \Omega$
  \begin{enumerate}
  \item[$(\Psi 1)'$]
    $\displaystyle
      (s_2-s_1)\bigl(\Psi(\,\cdot\,,s_2)-\Psi(\,\cdot\,,s_1)\bigr)
      \geq \sum_{i=1}^m \delta_i\,|s_2-s_1|^{r_i+1}\;.
    $
  \item[$(\Psi 2)'$]
    $\displaystyle
      N(s) - f\cdot 1_{\{\m(E)<\infty\}}
      \leq s\,\Psi(\,\cdot\,,s)
      \leq c \bigl( N(s)+ f\cdot 1_{\{\m(E)<\infty\}} \bigr)
    $ \\[+0.5ex]
    with $N$ as in \eqref{eq-3.10}.
  \item[$(\Phi 1)$]
    $\Phi(\,\cdot\,,s)=hs+\Phi_0(\,\cdot\,,s)$, $s\in\mathbb R$,
    for some $\scr F_t$-adapted
    $h\in L^\infty\bigl([0,T]\times \Omega;\mathrm dt\times P\bigr)$
    such that on $[0,T]\times \Omega$ we have
    \[
      \bigl| \Phi_0(\,\cdot\,,s_2) - \Phi_0(\,\cdot\,,s_1) \bigr|
      \leq \sum_{i=1}^m \delta_i\, \|L^{-1}\|_{\scr L(L^{r_i+1}(\m))}^{-1}\,
                       |s_2-s_1|^{r_i}\;.
    \]
  \item[$(\Phi 2)$]
    $\displaystyle
      \bigabs{\Phi_0(\argdot,s)}
      \leq \sum_{i=1}^m \tilde\varepsilon_i\,|s|^{r_i}\;,
    $\\[+0.5ex]
    where $\tilde \varepsilon_i
           := \varepsilon\,\|L^{-1}\|_{\scr L(L^{r_i+1}(\m))}^{-1}\, \varepsilon_i$
    for some $\varepsilon\in(0,1)$ (independent of $s\in \mathbb R$).
  \end{enumerate}
\end{enumerate}

\begin{rem}\label{rem-3.1}
  \begin{enumerate}
  \item[(i)]
    If \textbf{(A1)} holds for $\Psi$, then
    it holds for $s\mapsto -\Psi(-s)$ with the same Young
    function $N$.
  \item[(ii)]
    If $\m(E)=\infty$, then $(\Psi 2)$ and $(\Psi 3)$
    imply that $\Psi(\,\cdot\,,0)\equiv 0$.
  \item[(iii)]
    We note that if $(\scr E, D(\scr E))$ is a Dirichlet form
    and for some $c\in (0,\infty)$
    \[
      \scr E(u,u)
      \geq c\int u^2\d\m
      \quad\text{for all $u\in D(\scr E)$,}
    \]
    then as is well known, $(L, D(L))$ has bounded inverses as in
    condition \textbf{(A2)} for all $r_i\geq 1$.
  \end{enumerate}
\end{rem}

\begin{exa}\label{exa-3.5}
  \begin{enumerate}
  \item[(i)]
    Let $r_1,\dots,r_m > 0$, $\delta_1,\dots,\delta_m>0$
    and
    \[
      \Psi(s)
      := \operatorname{sign}(s)\sum_{i=1}^m \delta_i\, |s|^{r_i}\;.
    \]
    Define $N(s):= s\Psi(s)$. Then {\bf (A1)} is fulfilled,
     provided we can show that
    $N$, which is obviously a Young function, and also $N^*$
    are $\Delta_2$-regular.
    This is clear for $N$.
    To see the $\Delta_2$-regularity of $N^*$,
    let $r:=\min r_i$ and
    $\theta := 2^{1/r}>1$.
    We have, for $s>0$,
    \beg{equation}\label{eq-3.11}\beg{split}
      N^*(2s)&:=\sup_{t>0}\Big\{2st-\sum_{i=1}^m
     \dd_i t^{r_i+1}\Big\}=\sup_{t>0}\Big\{2\theta
      st-\sum_{i=1}^m\theta^{r_i+1}\dd_i t^{r_i+1}\Big\}\\
      &\le \theta^{r+1}\sup_{t>0}\Big\{2\theta^{-r}st-\sum_{i=1}^m
     \dd_i t^{r_i+1}\Big\}=\theta^{r+1}N^*(s).
    \end{split}\end{equation}
  \item[(ii)]
    Let $\Psi, N$ be as in (i), so that $(\Psi 2)'$ trivially holds, but with $r_1,\dots,r_m\geq 1$.
    Then an elementary calculation (see e.g.\
    \cite[p.~503]{Zei90}) shows that $(\Psi 1)'$ holds.
  \item[(iii)]
    Let $\Psi(s)
         := \operatorname{sign}(s)\,|s|^{\theta-1}
            \bigl( \log(1+|s|) \bigr)^r$, $s\in\mathbb R$,
    where $\theta\in(1,\infty)$, $r\in[1,\infty)$.
    Then $\Psi$ satisfies \textbf{\upshape (A1)} with
    $N(s) := s\,\Psi(s) = |s|^\theta \bigl(\log(1+|s|)\bigr)^r$,
    $s\in \mathbb R$.
    Obviously, $N$ is a Young function which together with its
    dual $N^*$ is easily checked to be $\Delta_2$-regular.
  \item[(iv)]
    Time dependent examples
    are easily obtained by e.g.\ multiplying $\Psi$ or $\Phi$ by a bounded
    adapted process which is bounded below by a strictly
    positive constant.
  \end{enumerate}
\end{exa}

To define $A$ as in \eqref{2.1}, we need one more lemma.

\begin{lem}\label{lem-3.6}
  \begin{enumerate}
  \item[(i)]
    Let \textbf{\upshape (A1)} hold.
    Then for all $s\in\mathbb R$
    \[
      N^*(c^{-1}\Psi(\,\cdot\,,s))
      \leq N(s) + \bigl[3 f + N^*\bigl(c^{-1}\Psi(\,\cdot\,,0)\bigr) \bigr]
                  \cdot 1_{\{\m(E)<\infty\}}
      \quad\text{on $[0,T]\times \Omega$.}
    \]
  \item[(ii)]
    Let \textbf{\upshape (A2)} hold and let $N$ be as in
    \eqref{eq-3.10}.
    Then there exists $\tilde c\in (0,\infty)$ such that for all
    $s\in\mathbb R$
    \[
      N^*\bigl(\Phi_0(\,\cdot\,,s) \bigr)
      \leq \tilde c \, N(s)
      \quad\text{on $[0,T]\times \Omega$.}
    \]
  \end{enumerate}
\end{lem}

\begin{proof} (i)    Let $s\in\mathbb R$.
    By Remark~\ref{rem-3.1}(i) we may assume that $s\geq 0$.
    Fix $t\in[0,T]$ and suppose first that
    $\Psi(t,s)\geq 0$,
    which by Remark~\ref{rem-3.1}(ii) and $(\Psi 1)$ is always the case
    if $\m(E)=\infty$.
    Then by $(\Psi 3)$
    \begin{align*}
      N^* (c^{-1} \Psi(t,s))
      &= \sup_{r\geq 0} \bigl[ rc^{-1}\Psi(t,s) -N(r) \bigr] \\
      &\leq \sup_{r\geq 0} \bigl[rc^{-1}(\Psi(t,s)-\Psi(\,\cdot\,,r))\bigr]
            + f\cdot  1_{\{\m(E)<\infty\}}  \\
      &= c^{-1} \sup_{s\geq r\geq 0} \bigl[r(\Psi(t,s)-\Psi(\,\cdot\,,r))\bigr]
            +f\cdot 1_{\{\m(E)<\infty\}} \;,
    \end{align*}
    where we used $(\Psi 1)$ in the last step.
    But since $N\geq 0$ the last term due to $(\Psi 2)$ is dominated
    by
    \[
      c^{-1} s\Psi(\,\cdot\,,s) + 2f \, 1_{\{\m(E)<\infty\}} \;.
    \]
    Now (i) follows by $(\Psi 3)$, for such $s\geq 0$, respectively
    is completely proved if $\m(E)=\infty$.

    If $\m(E)<\infty$ and $\Psi(t,s)<0$, then $|\Psi(t,s)|\leq |\Psi(t,0)|$
    by $(\Psi 1)$.
    Hence,
    \[
      N^*\bigl(c^{-1}\,\Psi(t,s)\bigr)
      =N^*\bigl(c^{-1} |\Psi(t,s)|\bigr)
      \leq N^*\bigl(c^{-1}|\Psi(t,0)|\bigr)
    \]
    and (i) also follows in this case.

  (ii)    Fix $1\leq i\leq m$ and define $N_i(s):= |s|^{r_i+1}$, $s\in\mathbb R$.
    Then
    \[
      N_i^*(s)
      = \biggl(\frac{1}{r_i+1}\biggr)^{\frac 1{r_i}}
        \frac{r_i}{r_i+1}\, |s|^{\frac{r_i+1}{r_i}}
      =: c_i\, |s|^{\frac{r_i+1}{r_i}}
    \]
    (as an elementary calculation shows).
    Hence, by $(\Phi 2)$ and \eqref{eq-3.11} for all $s\in\mathbb R$
    we have on $[0,T]\times \Omega$
    \begin{align*}
      N^*\bigl( \Phi_0(\,\cdot\,,s)\bigr)
      &= N^*\bigl(|\Phi_0(\,\cdot\,,s)|\bigr)
      \leq N^*\Bigl( \sum_{i=1}^m \tilde\varepsilon_i\,|s|^{r_i} \Bigr)
      \leq   N^* \Bigl(
                   \sum_{i=1}^m \tilde\varepsilon_i\,|s|^{r_i}
                 \Bigr) \\
      &\leq \sup_{\xi\geq 0}
              \Bigl[
                \xi\sum_{i=1}^m \tilde\varepsilon_i
                                \cdot \bigl(|s|^{r_i}-\xi^{r_i}\bigr)
              \Bigr]
      \leq \sum_{i=1}^m \tilde\varepsilon_i\,N_i^*\bigl(|s|^{r_i}\bigr)
      = \sum_{i=1}^m \tilde\varepsilon_i c_i\,|s|^{r_i+1}\;.
    \end{align*}
    Hence, the assertion follows with
    $\tilde c
     := \eps \max_{1\leq i\leq m}
             \bigl( c_i\,\norm{L^{-1}}_{\scr L(L^{r_i+1}(\m))}^{-1} \bigr)$.
    \end{proof}

Fix $\Psi,\Phi$ as above satisfying \textbf{(A1)} or \textbf{(A2)}.
We could define $A:[0,T]\times V\times\Omega\to V^*$ by
\begin{equation}\label{eq-3.12}
  A(t,v,\omega)
  := L\Psi(t,v,\omega)+\Phi(t,v,\omega)\;,
  \quad t\in[0,T],\; v\in V,\; \omega\in\Omega.
\end{equation}
By Lemma~\ref{lem-3.3}(iii) and Lemma~\ref{lem-3.6}(i)
the first summand is a well-defined element in $V^*$.
By Lemma~\ref{lem-3.6}(ii) this is also true for the second summand
since $V\subset V^*$ and $L_{N^*}=(L_N)^*\subset V^*$ since $V$
is dense in $L_N$ by \textbf{(N2)}.
But in order to show that $A$ satisfies our assumptions
\textbf{(H1)}--\textbf{(H4)}, we need estimates on $\Phi$
in \eqref{eq-3.12} and we have to compare it with $L\Psi$.
Therefore, we need to define the second summand in \eqref{eq-3.12}
in a more convenient way which is, however, only equivalent
under additional assumptions.
We refer to \cite{DRRW} for such a case and to the calculation
in Remark~\ref{rem-3.2} below. So, we define for
$t\in[0,T]$, $v\in V$, $\omega\in\Omega$
\begin{equation}\label{eq-3.13}
  \bar \Phi(t,v,\omega)
  := h_t\cdot v-\m\bigl( \Phi_0(t,v,\omega)\, L^{-1}\argdot \bigr)
\end{equation}
and
\begin{equation}\label{eq-3.14}
  A(t,v,\omega)
  := L\Psi(t,v,\omega) + \bar \Phi(t,v,\omega)\;.
\end{equation}
Since by assumption $\Phi_0$ is only nonzero if \textbf{(A2)}
holds, and then $L^{-1}:L_N\to L_N$ is continuous, it follows
by Lemma~\ref{lem-3.6}(ii) that
$V\ni u\mapsto \m\bigl(\Phi_0(t,v,\omega)\, L^{-1}u\bigr)$
is a continuous linear functional on $L_N$, hence on
$V=\scr F_e^*\cap L_N$, so belongs indeed to $V^*$.
With this definition of $A$ we shall then be able to verify
our conditions \textbf{(H1)}--\textbf{(H4)} on the basis
of our assumptions \textbf{(A1)} and \textbf{(A2)}.

\begin{rem}\label{rem-3.2}
  To avoid confusion below, we denote the continuous extension of
  the inverse $L^{-1}$ in \textbf{(A2)} to all of $L_N$ by
  $\tilde L^{-1}$.
  Now let $v\in V$ and $u\in V$ such that
  $\tilde L^{-1}u\in D(L)\cap L_{N^*}$.
  Then by Lemma~\ref{lem-3.3}(iii)
  \begin{align*}
    -\m(v\, \tilde L^{-1}u)
    &= \pairing{V^*}{\bar L \tilde L^{-1} u,v}{V}
    = \bracket{\bar L\tilde L^{-1}u,v}_{\scr F_e^*}
    = \bracket{u,v}_{\scr F_e^*}
    = \pairing{V^*}{v,u}{V}\;.
  \end{align*}
\end{rem}
So, if the set of such $u$ is dense in $V$, the above definitions
of $\Phi(t,v,\omega)$ and $\bar \Phi(t,v,\omega)$ are
equivalent (cf.\ \cite{DRRW} for an example), since $V$ is dense in
$L_N$.\\

Now we define
\begin{equation}\label{eq-3.15}
  R(v)
  := \m(N(v)) + \norm v_H^2\;,\quad v\in V,
\end{equation}
and
\begin{equation}\label{eq-3.16}
  K
  := L_N \bigl( [0,T]\times E\times\Omega;\mathrm dt\times \m\times P \bigr)
     \cap L^2\bigl( [0,T]\times\Omega\to H; \mathrm dt\times P \bigr)
\end{equation}
with norm
\begin{equation}\label{eq-3.17}
  \norm{\argdot}_K
  := \norm{\argdot}_{L_N ( [0,T]\times E\times\Omega;\mathrm dt\times \m\times P ) }
     + \norm{\argdot}_{L^2 ( [0,T]\times\Omega\to H; \mathrm dt\times P ) }\;.
\end{equation}
That the intersection in \eqref{eq-3.16} is meaningful follows from
the last inclusion in the following lemma and the definition of $V$
($= L_N(\m)\cap \scr F_e^*$). It also follows that $K$ is complete
(hence a Banach space) and reflexive.

\begin{lem}\label{lem-3.7}
  Let $N$ be a $\Delta_2$-regular Young function and $q\in(2,\infty)$
  as in Lemma~\ref{lem-3.2}.
  Then the following embeddings are dense and continuous
  \begin{align*}
    &L^q\bigl([0,T]\times\Omega\to V;\mathrm dt\times P\bigr)
    \subset L^q\bigl([0,T]\times\Omega\to L_N(\m);\mathrm dt\times P\bigr)\\
    &\quad\subset L_N \bigl(
                        [0,T]\times E\times \Omega;
                        \mathrm dt\times \m\times P
                      \bigr)
    \subset L^1\bigl([0,T]\times \Omega\to L_N(\m);\mathrm dt\times P\bigr)\;.
  \end{align*}
\end{lem}

\begin{proof}
  The assertion with respect to the first inclusion is clear, by
  (the second half of) condition \textbf{(N2)}.
  To prove the assertion for the second inclusion, let
  $g\in L_{N^*}\bigl([0,T]\times E\times \Omega;\mathrm dt\times \m\times P\bigr)$
  such that
  \[
    \bar \m (N(g))
    \quad\biggl(:= \int N(g)\;\mathrm d\bar \m\biggr)
    \quad \leq 1\;,
  \]
  where $\bar \m := \mathrm dt\times\m\times P$.
  Then for all $f\in L_N(\bar\m)$, since $s_1\cdot s_2\leq N(s_1)+N^*(s_2)$
  for all $s_1,s_2\in \R$,
  \begin{align*}
    \bar\m(f\cdot g)
    &\leq \bar \m\bigl(N(f)\bigr) + 1\\
    &\leq \int\!\int_0^T
          \m\Biggl(
           N \biggl(
               \bigl( \norm{f(t,\argdot,\omega)}_{L_N(\m)} +2 \bigr)
               \,\frac{f(t,\argdot,\omega)}
                      {\norm{f(t,\argdot,\omega)}_{L_N(\m)} +2}
             \biggr) \Biggr)
           \;\mathrm dt\; P(\mathrm d\omega)
         + 1\\
    &\leq \int\!\int_0^T
            \bigl(\norm{f(t,\argdot,\omega)}_{L_N(\m)}+2\bigr)^q
            \;\bigl( 1+2\m(E)\, 1_{\{\m(E)<\infty\}} \bigr)
            \;\mathrm dt \;P(\mathrm d\omega)
          + 1\;,
  \end{align*}
  where we used \eqref{eq:delta-reg} in the last step.
  Now by \eqref{*0} we obtain that for some constants $a,b>0$
  (independent of $f$)
  \[
    \norm f_{L_N(\bar\m)}
    \leq a \int\!\int_0^T \bignorm{f(t,\argdot,\omega)}_{L_N(\m)}^q
                         \;\mathrm dt\; P(\mathrm d\omega)
         + b \;.
  \]
  Hence, $\norm{\argdot}_{L_N(\bar\m)}$ is bounded on bounded sets
  of $L^q\bigl( [0,T]\times\Omega\to L_N(\m);\mathrm dt\times P \bigr)$,
  so the second embedding in the assertion is  continuous.
  To show its density, it is enough to prove that
  \[
    L_0^\infty(\bar \m)
    \subset L^q\bigl([0,T]\times\Omega\to L_N(\m);\mathrm dt\times P\bigr)\;,
  \]
  where $L_0^\infty(\bar\m)$ denotes the set of all $f\in L^\infty(\bar \m)$
  such that for some $E_0\in\scr B$ with $\m(E_0)<\infty$,
  $\{f\not= 0\}\subset [0,T]\times E_0\times\Omega$.
  (Obviously, $L_0^\infty(\bar\m)$ is dense in $L_N(\bar \m)$.)
  But by Lemma~\ref{lem-3.2}(ii) there exists $c\in(0,\infty)$
  such that for all $f\in L_0^\infty(\bar \m)$
  \begin{align*}
    &\int\!\int \bignorm{f(t,\argdot,\omega)}_{L_N(\m)}^q
      \;\mathrm dt\;P(\mathrm d\omega)\\
    &\quad\leq c \int\!\int
                   \bigl(
                     \norm{f(t,\argdot,\omega)}_{L^q(\m)}^q
                     + \norm{f(t,\argdot,\omega)}_{L^1(\m)}^q
                   \bigr)
                   \;\mathrm dt\; P(\mathrm d\omega)\\
    &\quad\leq c\,\norm f_{L^\infty(\bar \m)}^q\;
               T\bigl( \m(E_0)+\m(E_0)^q \bigr)
        <\infty \;.
  \end{align*}

  Now let us prove the assertion for the last inclusion.
  Let $f\in L_N(\bar \m)\setminus\{0\}$.
  As above by \eqref{*0} we obtain that for $\mathrm dt\times P$-a.e.\
  $(t,\omega)\in[0,T] \times\Omega$
  \[
    \bignorm{f(t,\argdot,\omega)}_{L_N(\m)}
    \leq \m\bigl( N(f(t,\argdot,\omega)) \bigr) +1\;.
  \]
  Hence,
  \begin{align*}
    &\int\!\int \bignorm{f(t,\argdot,\omega)}_{L_N(\m)}\;\mathrm dt\;P(\d\omega)
    \leq \bar\m(N(f)) +T \\
    &\quad\leq \bar\m \Biggl( N \biggl(
                        \frac{f}{\norm f_{L_N(\bar \m)}}
                        \; \bigl[\norm f_{L_N(\bar \m)}+2\bigr]
                      \biggr)\Biggr)
               +T \\
    &\quad\leq \bigl(\norm f_{L_N(\bar \m)}+2\bigr)^q
         \, \bigl( 1+2T\,\m(E)\, 1_{\{\m(E)<\infty\}} \bigr) +T\;,
  \end{align*}
  where we used \eqref{eq:delta-reg} in the last step.
  Since the last expression is bounded on bounded sets of $L_N(\bar \m)$,
  the third continuous embedding in the assertion is proved.
  Its density is then obvious since clearly
  $L_0^\infty(\bar\m)\subset L^1\bigl([0,T]\times \Omega\to L_N(\m);
                                     \mathrm dt\times P\bigr)$.
\end{proof}

By definition and Lemma~\ref{lem-3.7} it now follows that for $q$
as in Lemma~\ref{lem-3.2}
\[
  L^q\bigl( [0,T]\times \Omega\to V;\mathrm dt\times P \bigr)
  \subset K \subset
  L^1\bigl( [0,T]\times\Omega\to V;\d t\times P \bigr)
\]
continuously, and both embeddings are dense, since $L^q$ is dense
in $L^1$.
So, it remains to check our general conditions \textbf{(K)},
\textbf{(H1)}--\textbf{(H4)}.

\begin{prp}\label{prop-3.8}
  For $R$ and $K$ defined in \eqref{eq-3.15}, \eqref{eq-3.16},
  respectively, conditions \textbf{\upshape (K)},
  \textbf{\upshape (H1)}--\textbf{\upshape (H4)}
  from Section~\ref{sec:mon-sdes} hold for $A$ defined in
  \eqref{eq-3.14} and \eqref{eq-3.13}.
\end{prp}

\begin{proof}
  To prove \textbf{(K)}(i) it suffices to show that for any
  sequence $z^{(n)}\in K$, $n\in\N$, one has
  $\norm{z^{(n)}}_{L_N(\bar\m)}\to 0$ if and only if
  $\bar\m(N(z^{(n)}))\to 0$, where as before
  $\bar \m := \d t\times\m \times P$.
  So, assume $\norm{z^{(n)}}_{L_N(\bar\m)}\to 0$.
  Then $z^{(n)}\to 0$ in $\bar\m$-measure, because for all $k\in\N$
  there exists $n_k\in \N$ such that
  \[
    \bar\m \Biggl( N \biggl( \frac{\abs{z^{(n_k)}}}{2^{-k}} \biggr)\Biggr)
    \leq 1\;,
  \]
  hence by the convexity and continuity of $N$
  \[
    \bar \m\biggl( N\Bigl( \sum_{k=1}^\infty \abs{z^{(n_k)}} \Bigr) \biggr)
    \leq \sum_{k=1}^\infty 2^{-k}\,\bar\m\bigl(N(2^k\,\abs{z^{(n_k)}})\bigr)
    < \infty\;,
  \]
  so $\sum_{k=1}^\infty \abs{z^{(n_k)}}<\infty$ $\bar\m$-a.e.
  Since this is true for any subsequence of $z^{(n)}$, $n\in\N$,
  this really implies that $z^{(n)}\to 0$ in $\bar\m$-measure.

  On the other hand, by \eqref{eq:delta-reg} for
  $N^*$ (with possibly different $q>2$), for any $\vv\in (0,
  2^{1-q}]$ we have
  \begin{equation*}
    \begin{split}
      N(\vv s)
      &:= \sup_{r>0} \{|s|\vv r -N^*(r)\}
      = \vv^{q/(q-1)}\sup_{r>0} \{|s|\vv^{-1/(q-1)}r- N^*(r) \vv^{-q/(q-1)}\}\\
      &\le \vv^{q/(q-1)}\sup_{r>0} \{|s|\vv^{-1/(q-1)} r - N^*(\vv^{-1/(q-1)} r)\}
           + 2\cdot 1_{\{\m(E)<\infty\}}\\
      &= \vv^{q/(q-1)} N(s) + 2\cdot 1_{\{\m(E)<\infty\}} ,\ \ \ s\in\R.
    \end{split}
  \end{equation*}
  This implies, for ($0\not=$) $\|z^{(n)}\|_{L_N(\bar\m)}\le 2^{1-q},$
  \[
    N(z^{(n)})
    = N \biggl(
          \ff{z^{(n)} \|z^{(n)}\|_{L_N(\bar\m)}}
             {\|z^{(n)}\|_{L_N(\bar\m)}}
        \biggr)
    \le \|z^{(n)}\|_{L_N(\bar\m)}^{q/(q-1)} \;
        N \biggl( \ff{z^{(n)}}{\|z^{(n)}\|_{L_N(\bar\m)}} \biggr)
        + 2T\,\m(E)\, 1_{\{\m(E)<\infty\}}\;.
  \]
  But the right hand side converges in $L^1(\bar\m)$, hence $N(z^{(n)})\to 0$
  in $L^1(\bar\m)$.

  Now assume that $N(z^{(n)})\to 0$ in $L^1(\bar\m)$, hence in
  $\bar\m$-measure.
  Then $z^{(n)}\to 0$ in $\bar\m$-measure.
  But for $\lambda\in\bigl(0,\frac 12\bigr)$ by \eqref{eq:delta-reg}
  we have
  \[
    N\biggl(\frac{z^{(n)}}{\lambda}\biggr)
    \leq \lambda^{-q} \bigl(N(z^{(n)}) + 2\cdot 1_{\{\m(E)<\infty\}}\bigr)
  \]
  and the right hand side converges in $L^1(\bar\m)$.
  Hence, $N\bigl(\frac{z^{(n)}}{\lambda}\bigr)\to 0$ in $L^1(\m)$, so
  for sufficiently large $n$
  \[
    \norm{z^{(n)}}_{L_N(\bar\m)} \leq \lambda\;,
  \]
  thus $\norm{z^{(n)}}_{L_N(\bar\m)}\to 0$ since
  $\lambda\in\bigl(0,\frac 12\bigr)$ was arbitrary.

  Now we verify \textbf{(K)}(ii).
  By \eqref{eq:delta-reg} and since $N(s)$ is increasing in $\abs s$,
  we have for $z\in K\setminus \{0\}$
  \begin{align*}
    \bar\m\bigl(N(z)\bigr)
    &\leq \bar\m \Biggl( N \biggl(
                   \frac{z(\norm z_{L_N(\bar\m)} +2)}
                        {\norm z_{L_N(\bar\m)}}
                 \biggr)\Biggr)\\
    &\leq \bar\m
            \Biggl(
              \biggl[
                N \biggl(
                    \frac{z}{\norm z_{L_N(\bar\m)}}
                  \biggr)
                + 2\cdot 1_{\{\m(E)<\infty\}}
              \biggr]\;
              \bigl(
                \norm z_{L_N(\bar\m)}+ 2
              \bigr)^q
            \Biggr) \\
    &\leq \bigl(
            1+2T\,\m(E)\, 1_{\{\m(E)<\infty\}}
          \bigr)\,
          \bigl(
            \norm z_{L_N(\bar\m)} +2
          \bigr)^q \;.
  \end{align*}
  This implies that for some $c\in(0,\infty)$
  \begin{align*}
    \norm z_K
    &\geq c \Biggl(
              \bigl[
                \bar\m\bigl(N(z)\bigr)
              \bigr]^{1/q}
              + \biggl(
                  \E \int_0^T \norm{z_t}_H^2\;\mathrm dt
                \biggr)^{1/2}
            \Biggr) -2\\
    &\geq c \biggl(
              \E\int_0^T R(z_t)\;\d t
            \biggr)^{1/q}
          - c - 2\;,
  \end{align*}
  where we used the elementary estimate
  $(a+b)^{1/q}\leq a^{1/q} + b^{1/2}+1$,
  $a,b\geq 0$, $q\geq 2$.
  On the other hand, since by \eqref{*0} and because of
  $st\leq N(s)+N^*(t)$ we have
  \[
    \norm z_{L_N(\bar\m)}
    \leq \bar\m\bigl(N(z)\bigr) +1\;,
  \]
  it follows that
  \[
    \norm z_K
    \leq \E\int_0^T \m\bigl(N(z_t)\bigr)\;\d t
         + 1 + \biggl( \E\int_0^T\norm{z_t}_H^2\;\d t \biggr)^{1/2}
    \leq \E \int_0^T R(z_t)\;\d t + 2\;.
  \]
  Therefore, \textbf{(K)}(ii) holds for
  \[
    W_1(r)
    := cr^{1/q} -c-2\;,\quad
    W_2(r) := r+2\;,\quad r\geq 0.
  \]
 Since $N$ is
convex,

\beg{equation*}\beg{split}R(x+y)&:= \mathbf m
(N(x+y))+\|x+y\|_H^2\\
&\le \ff 1 2 \mathbf m(N(2x)+N(2y)) + 2 \|x\|_H^2+2\|y\|_H^2=\ff 1 2
(R(2x)+R(2y)).\end{split}\end{equation*} So, {\bf (K)}(iii) holds
for $C=\ff 1 2.$ \textbf{(K)}(iv) is clear, since obviously
  for $z\in L_N(\bar\m)$ and
  $h\in L^\infty\bigl([0,T]\times\Omega;\mathrm dt\times P\bigr)$
  \[
    \norm{hz}_{L_N(\bar\m)}
    \leq \norm h_{L^\infty(\mathrm dt\times P)}\,\norm z_{L_N(\bar \m)}\;.
  \]

  Now we are going to prove \textbf{(H1)}--\textbf{(H4)}.
  \begin{enumerate}
  \item[\textbf{(H1)}:]
    Let $u,v,x\in V=L_N(\m)\cap\scr F^*_e$ and $\lambda\in\R$.
    Then by \eqref{eq-3.11}, \eqref{eq-3.10}, Lemma~\ref{lem-3.6}
    and \eqref{eq-3.7} on $[0,T]\times \Omega$ we have
    \begin{align*}
      &\bigpairing{V^*}{A(\argdot,u+\lambda v),x}{V}\\
      &\quad= - \m\bigl( \Psi(\argdot,u+\lambda v) \,x \bigr)
        + h_t \pairing{V^*}{u+\lambda v,x}{V}
        - \m\bigl( \Phi_0(\argdot,u+\lambda v) \, L^{-1}x \bigr) \;.
    \end{align*}
    But by Lemma~\ref{lem-3.6} on $[0,T]\times \Omega$
    \begin{align*}
      \bigabs{\Psi(\argdot,u+\lambda v)\cdot x}
      &\leq N^*\bigl( c^{-1}\, \Psi(\argdot,u+\lambda v) \bigr) +N(cx)\\
      &\leq N(u+\lambda v)
           + \bigl[ 3f + N^*\bigl(c^{-1}\, \Psi(\argdot , 0)\bigr) \bigr]
             \cdot 1_{\{\m(E)<\infty\}}
           + N(cx)
    \end{align*}
    and
    \begin{align*}
      \bigabs{\Phi_0(\argdot,u+\lambda v)}\,\abs{L^{-1}x}
      \leq N^*\bigl(\Phi_0(\argdot,u+\lambda v)\bigr) + N(L^{-1} x)
      \leq \tilde c\, N(u+\lambda v) + N(L^{-1}x)\;,
    \end{align*}
    where for $\lambda \in[-1,1]$
    \[
      N(u+\lambda v)
      = N\bigl(\abs{u+\lambda v}\bigr)
      \leq N\bigl( \abs u + \abs v\bigr)\;.
    \]
    So \textbf{(H1)} follows by the continuity of $\Psi,\Phi_0$ in
    the spatial variable and Lebesgue's dominated convergence theorem.
  \item[\textbf{(H2)}:]
    Let $u,v\in V$. Then as above on $[0,T]\times \Omega$
    \begin{equation}\label{eq-3.18}
      \!\!\begin{aligned}
        &\bigpairing{V^*}{A(\argdot,u)-A(\argdot,v), u-v}{V}\\
        &\quad= -\m \Bigl(
                      \bigl( \Psi(\argdot,u)-\Psi(\argdot,v) \bigr)
                      \, (u-v)
                    \Bigr)
                {}+ h_t\,\pairing{V^*}{u-v,u-v}{V} \\
        &\qquad{}-\m \Bigl(
                      \bigl( \Phi_0(\argdot,u)-\Phi_0(\argdot,v) \bigr)
                      \, L^{-1}(u-v)
                    \Bigr) \;.
      \end{aligned}
    \end{equation}
    In case \textbf{(A1)} holds, the latter is dominated by
    $\norm h_{L^\infty(\mathrm dt\times P)}\norm{u-v}_H^2$.
    In case \textbf{(A2)} holds, the absolute value of the last
    summand is by $(\Phi 1)$ dominated by
    \begin{align*}
      &\sum_{i=1}^m \delta_i \norm{L^{-1}}_{\scr L(L^{r_i+1}(\m))}^{-1}
        \,\m \bigl(
               \abs{u-v}^{r_i}
               \,\bigabs{L^{-1}(u-v)}
             \bigr)\\
      &\quad\leq \sum_{i=1}^m \delta_i
                    \Bigl(\m(\abs{u-v}^{r_i+1})\Bigr)^{r_i/(r_i+1)}
         \, \m(\abs{u-v}^{r_i+1})^{1/(r_i+1)}\\
      &\quad\leq \m \Bigl(
                \bigl(\Psi(\argdot,u)-\Psi(\argdot,v)\bigr)\,(u-v)
              \Bigr)\;,
    \end{align*}
    where we first used H\"older's inequality and then $(\Psi 1)'$.
    So, also in case \textbf{(A2)} holds, the right hand side of
    \eqref{eq-3.18} is dominated by
    $\norm h_{L^\infty(\mathrm dt\times P)}\norm{u-v}_H^2$ on $[0,T]\times\Omega$,
    so \textbf{(H2)} is proved.
  \item[\textbf{(H3)}:]
    Let $v\in V$. Then as above on $[0,T]\times\Omega$,
    \[
      \bigpairing{V^*}{A(\argdot,v),v}{V}
      = - \m \bigl(\Psi(\argdot,v)\,v\bigr)
        + h_t\,\pairing{V^*}{v,v}{V}
        - \m\bigl( \Phi_0(\argdot,v)\, L^{-1}v \bigr)\;.
    \]
    By $(\Psi 2)$, $(\Psi 2)'$, respectively, and $(\Phi 2)$, this is on
    $[0,T]\times \Omega$ dominated by
    \[
      - (1-\eps)\,\m\bigl(N(v)\bigr)
      + f\, \m(E)\, 1_{\{\m(E)<\eps\}}
      + \norm h_{L^\infty(\mathrm dt\times P)}\,\norm v_H^2\;.
    \]
    By the definition of $R$ (cf.\ \eqref{eq-3.15}) condition \textbf{(H3)}
    now follows.
  \item[\textbf{(H4)}:]
    Let $u,v\in V$. Then as above on $[0,T\times \Omega]$
    \begin{align*}
      &\bigabs{\pairing{V^*}{A(\argdot,v),u}{V}}\\
      &\quad\leq \m\bigl(\abs{\Psi(\argdot,v)}\,\abs u\bigr)
        + \norm h_{L^\infty(\mathrm dt\times P)}\,\norm v_H\,\norm u_H
        + \m\bigl( \abs{\Phi_0(\argdot,v)}\,\abs{L^{-1}u} \bigr)\;.
    \end{align*}
    But by Lemma~\ref{lem-3.6} on $[0,T]\times\Omega$
    \begin{align*}
      c^{-1}\,\bigabs{\Psi(\argdot,v)}\,\abs u
      &\leq N^*\bigl(c^{-1}\,\Psi(\argdot,v)\bigr) +N(u)\\
      &\leq N(v)+ \bigl[ 3f+N^*\bigl(c^{-1}\,\Psi(\argdot,0)\bigr) \bigr]
                  \cdot 1_{\{\m(E)<\infty\}}
            + N(u)\;.
    \end{align*}
    Furthermore, by $(\Phi 2)$, H\"older's and Young's inequality
    \begin{align*}
      \m \bigl( \abs{\Phi_0(\argdot,v)}\, \abs{L^{-1}u} \bigr)
      &\leq \eps\sum_{i=1}^m\eps_i\bigl(\m(\abs v^{r_i+1})\bigr)^{r_i/(r_i+1)}
        \,\bigl( \m(\abs u^{r_i+1}) \bigr)^{1/(r_i+1)} \\
      &\leq \eps\sum_{i=1}^m
                  \eps_i\biggl(
                          \frac{r_i}{r_i+1}\;\m\bigl(\abs v^{r_i+1}\bigr)
                          + \frac 1{r_i+1}\;\m\bigl(\abs u^{r_i+1}\bigr)
                        \biggr)\\
      &\leq \eps\bigl(N(v)+N(u)\bigr)\;,
    \end{align*}
    and \textbf{(H4)} follows.\qedhere
  \end{enumerate}
\end{proof}
Now our general results from Section~\ref{sec:mon-sdes} apply to
this case.

\begin{thm}\label{thm-summary}
  Assume that conditions \textbf{\upshape (L1)},
  \textbf{\upshape (N1)} and \textbf{\upshape (N2)}
  hold and let $\Phi,\Psi$ be as in \eqref{eq-3.9a} satisfying
  \textbf{\upshape (A1)} or \textbf{\upshape (A2)}
  for some Young function $N$ and let $B$ be as in \eqref{eq-3.9b}.
  Then for any $X_0\in L^2(\Omega\to H,\scr F_0; P)$ the SPDE
  \[
    \mathrm d X_t
    = \bigl[
        L\Psi(t,X_t)
        + \bar \Phi (t,X_t)
      \bigr] \;\mathrm d t
    + B (t,X_t) \;\mathrm d W_t
  \]
  has a unique continuous $H$-valued solution (in the sense of
  Definition~2.1 with $H$ being the Green space of $L$ and
  $K$ defined by \eqref{eq-3.16}) with initial condition $X_0$.
  The solution satisfies
  \[
    \E \sup_{t\in[0,T]} \norm{X_t}_H^2
    <\infty
  \]
  and all further assertions from Theorem~\ref{T2.1}
  and Proposition~\ref{P2.2} also hold.
\end{thm}

\section{Appendix: the It\^o formula for the square of the norm}

In this section we aim to prove the It\^o formula for
$\|X_t\|_H^2$ which has been used in the paper. This formula has
been established in \cite{KR} for the case where $K:=
L^p([0,T]\times \OO\to V; \d t\times P)$ for some $p>1,
R:=\|\cdot\|_V^p, W_1=W_2:=\|\cdot\|^{1/p}.$

As in \cite{KR}, let us first consider the piecewise-constant
approximation of a process in $K$ by using an argument of Doob
(see \cite[Ch. IX, \S 5]{Doob}, or \cite[page 75]{R}).

\begin{lem}\label{A1}
Assume {\bf (K)} and let $X: [0,T]\times
\OO\to V^*$ be $\scr B([0,T])\times \scr F/ \scr B(V^*)$-measurable
such that $X=\bar X\ \d t\times P$-a.e. for some
$\mathrm dt\times P$-version $\bar X$ of an element in $K$.
Then there exists a sequence of  partitions
$I_l:=\{0=t_0^l<t_1^l<\cdots<t_{k_l}^l=T\}$ such that $I_l\subset
I_{l+1}$ and  $\dd(I_l):=\max_i (t_i^l-t_{i-1}^l)\to 0$ as
$l\to\infty$ and for

$$\bar X^l:= \sum_{i=2}^{k_l}1_{[t_{i-1}^l, t_i^l)}X_{t_{i-1}^l},
\ \ \tt X^l:= \sum_{i=1}^{k_l-1}1_{[t_{i-1}^l, t_i^l)}X_{t_{i}^l},
\ \ \ l\ge 1$$ we have $\bar X^l,\tt X^l$ are $V$-valued $P$-a.s.
and $(\d t\times P$-versions of$)$ elements in $K$. Furthermore,

$$\lim_{l\to\infty}\big\{\|\bar X-\bar X^l\|_K+\|\bar X-\tt
X^l\|_K\big\}=0.$$ In particular, $X_{t_i^l}\in V$ for all $l\ge
1, 1\le i\le k_l-1.$
\end{lem}

\begin{proof}
For simplicity we  assume that $T=1$ and let $X$ be
extended to $\R\times\Omega$ by setting $X|_{[0,1]^c}=0.$ Since
$L^p([0,T]\times \OO\to V; \d t\times P)$ is dense in $K$, it is
easy to see from {\bf (K)} that there exists $\OO'\subset \OO$
with full probability such that $\int_0^1R(\bar X_t)\d t)<\infty$
holds on $\OO'$ and for any $\oo\in \OO',$ there exists a sequence
$\{f_n\}\subset C(\R;V)$ with compact support such that
$$\int_\R R(4f_n(t)- 4\bar X_t(\oo))\d t\le \ff 1 n,\ \ \ \ n\ge 1.$$
Thus for every $n\ge 1$ it follows from {\bf (K)}(iii), since
$R(v)\to 0$ as $v\to 0$, and $f_n\in C_0(\R;V)$, that

\beg{equation*}\beg{split} &\limsup_{s\to\infty} \int_R R(\bar
X_{t+1}(\oo)-\bar X_t(\oo))\d t\\
&\le C\limsup_{s\to\infty} \int_R\big\{R(2\bar
X_{s+t}(\oo)-2f_n(t+s)+ 2\bar X_t(\oo) -2 f_n(t)) + R(2 f_n(t+s) -2
f_n(t))\big\}\d t\\
& \le C^2 \limsup_{s\to\infty} \int_R \big\{R(4\bar
X_{s+t}(\oo)-4f_n(t+s))+ R(4\bar X_t(\oo) -4f_n(t))\big\}\d t\le
\ff{2C^2}n.\end{split}\end{equation*} Note here that since by
continuity $R$ is bounded on sufficiently small balls around $0$ and
since each $f_n$ is uniformly continuous we really have by dominated
convergence that for all $n\ge 1$

$$\lim_{s\to 0} \int_\R R(2f_n(s+t)-2f_n(t))\d t =0.$$
Letting $n\to\infty$ we arrive at

\beq\label{A.1}\lim_{\dd\to 0}\int_\R R(\bar X_{\dd+s}(\oo)-\bar
X_s(\oo))\d s =0,\ \ \ \oo\in \OO'.\end{equation} Now, given $t\in
\R$, let $[t]$ denote the biggest integer $\le t$. Let $\gg_n(t):=
2^{-n}[2^nt],\ n\ge 1.$ Shifting the integral in (\ref{A.1}) by
$t$ and taking $\dd=\gg_n(t)-t$ we obtain

$$\lim_{n\to \infty} \int_\R R(\bar X_{\gg_n(t)+s}-\bar X_{t+s})\d s=0\ \
\text{on}\ \OO'.$$ Moreover,  since $R(0)=0$ and by {\bf (K)}(iii)
and Remark 2.1(3)

\beg{equation*}\beg{split}\int_0^1 R(\bar X_{\gg_n(t)+s}-\bar
X_{t+s})\d s&\le 1_{[-2,2]}(t)C\int_\R \big\{R(2\bar
X_{\gg_n(t)+s})+R(2\bar X_{t+s})\big\}\d
s\\
&= 2C1_{[-2,2]}(t)\int_0^1 R(2\bar X_s)\d
s<\infty.\end{split}\end{equation*} So, by the dominated convergence
theorem, we obtain that

\beq\label{A.2}0= \lim_{n\to\infty}\E\int_\R \d t \int_0^1 R(\bar
X_{\gg_n(t)+s}-\bar X_{t+s})\d s \ge\lim_{n\to\infty}\E\int_0^1\d
s \int_0^1 R(\bar X_{\gg_n(t-s)+s}-\bar X_t)\d t\;.
\end{equation}
Given $s\in [0,1)$ and $n\ge 1$, let the
partition $I_n(s)$ be defined by

$$t_0^n(s):=0,\ t_i^n(s):= \Big(s- \ff{[2^n s]}{2^n}\Big) + \ff{i-1}{2^n},\ 1\le i\le 2^n,\
t_{2^n+1}^n(s):=1.$$ Then, for $t\in [t_{i-1}^n(s), t_i^n(s))$ one
has $t-s\in [2^{-n}(i-[2^n s]-2), 2^{-n}(i-[2^n s]-1))$ and hence,

$$\gg_n(t-s)+s=\big\{2^{-n}(i-[2^n s]-2)+s\big\}^+ = t_{i-1}^n(s),\ \  1\le i\le 2^n +2.$$
Therefore, (\ref{A.2}) implies

$$\lim_{n\to\infty}\E\int_0^1\d s\int_0^1 R(\bar X_t-\bar X_t^{n,s})\d t
=0,$$ where $\bar X^{n,s}$ is the process defined as $\bar X^l$ for
the partition $I_n(s)$ but with $X_{t_i^l(s)}$ replaced by $\bar
X_{t_i^l(s)}$. Similarly, the same holds for $\tt X^{n,s}$ in place
of $\bar X^{n,s}$ by using $\tt\gg_n:= \gg_n+2^{-n}$ instead of
$\gg_n$, where $\tt X^{n,s}$ is defined as $\tt X^l$ for the
partition $I_n(s)$ but with $X_{t_i^l(s)}$ replaced by $\bar
X_{t_i^l(s)}$. Hence, there exist a subsequence $n_k\to\infty$ and a
$\d s$-zero set $N_1\in \scr B([0,1])$ such that

$$\lim_{k\to\infty} \E\int_0^1 \big\{R(\bar X_t-\bar X_t^{n_k,s})+
R(\bar X_t- \tt X_t^{n_k,s})\big\}\d t=0,\ \ \ s\in [0,1]\setminus
N_1.$$ Since for $1\le i\le 2^n$ the maps $s\mapsto t_i^n(s)$ are
piecewise $C^1$-diffeomorphisms, the image measures of $\d s$
under these maps are absolutely continuous with respect to $\d s$.
Therefore, since $\bar X= X\ \d s\times P$-a.e., there exists a
$\d s$-zero set $N_2\in \scr B([0,1])$ such that

$$\bar X_{t_i^n(s)} = X_{t_i^n(s)}\ \ \text{a.s.},\ \ \
s\in [0,1]\setminus N_2, 1\le i\le 2^n, n\geq 1.$$ Since for any $s\in
[0,1]\setminus (N_1\cup N_2)$ one has $\E R(\bar
X_{t_i^n(s)})<\infty$ and by Remark 2.1(3) $z\in K$ if and only if
$z\in L^1([0,1]\times \OO; \d t\times P)$ with $\E\int_0^1
R(z_t)\d t<\infty,$ the map

$$[0,1]\times \OO\ni (t,\oo)\mapsto z^{i,n}_t:=X_{t_i^n(s)}(\oo)\in V$$
is once again in $K$. Therefore, fixing $s\in [0,1]\setminus
(N_1\cup N_2)$, the sequence of the corresponding partitions
$I_{n_l}(s), l\ge 1$, has all properties of the assertion.
\end{proof}

\paragraph{Remark 4.1.} As follows from the above proof all the
partition points $t_i^l, l\ge 1, 1\le i\le k_l-1,$ in the
assertion of Lemma \ref{A1} can be chosen outside an a priori
given Lebesgue zero set in $[0,T].$

\beg{thm}\label{TA}
Assume {\bf (K)}. Let $X_0\in L^2(\OO\to H; \F_0; P)$,  $Y\in
K^*$, $Y$ progressively measurable, and $Z\in J$. Define the
continuous $V^*$-valued adapted process

$$X_t:=X_0 +\int_0^t Y_s \d s +\int_0^t Z_s\d W_s,\ \ \ t\in
[0,T].$$
If there exists a $\mathrm dt\times P$-version $\bar X$ of an
element in $K$ such that $X=\bar X\ \d
t\times P$-a.e., then $X_t$ is a
continuous process on $H$ such that $\E\sup_{t\le
T}\|X_t\|_H^2<\infty$ and $P$-a.s.

\beq\label{Ito}\|X_t\|_H^2 =\|X_0\|_H^2 +\int_0^t (2\
_{V^*}\<Y_s,\bar X_s\>_V +\|Z_s\|_{\L_{HS}}^2)\d s + 2 \int_0^t
\<Z_s\d W_s, X_s\>_H,\ \ t\in [0,T].
\end{equation}
Setting $_{V^*}\<Y_s, X_s\>_V=0$ for $X_s\notin V$
we may replace in the right-hand side $\bar X_s$ by $X_s.$
\end{thm}

\begin{proof}
Since  $M_t:=\int_0^t Z_s\d W_s$ is already a
continuous martingale on $H$ and since $Y\in K^*\subset
L^{p/(p-1)}([0,T]\times\OO\to V^*;\d t\times P)$ is progressively
measurable, $\int_0^t Y_s\d s$ is a continuous adapted process on
$V^*$. Thus, $X$ is a continuous adapted process on $V^*$, hence  is
$\scr B([0,T])\times\scr F/\scr B(V^*)$-measurable.
Then, due to Lemma \ref{A1} and Remark~2.1(5), the
remainder of the proof is similar to that of \cite[Theorem
I.3.1]{KR}. We include a complete proof below  for the readers'
convenience.

(a) Note that by \cite[Lemma I.4.2]{KR},

\beq\label{Tr1}\beg{split} \|X_t\|_H^2=&\|X_s\|_H^2 + 2 \int_s^t\
_{V^*}\<Y_r, X_t\>_V\d r + 2 \<X_s, M_t-M_s\>_H\\
& +\|M_t-M_s\|_H^2-\|X_t-X_s-M_t+
M_s\|_H^2\end{split}\end{equation} holds for all $t>s$ such that
$X_t,X_s\in V.$ Indeed, this follows immediately by noting that

\beg{equation*}\beg{split}
&\|M_t-M_s\|_H^2-\|X_t-X_s-M_t+M_s\|_H^2 +2\<X_s, M_t-M_s\>_H\\
&= 2\<X_t, M_t-M_s\>_H -\|X_t-X_s\|_H^2\\
&= 2\<X_t, X_t-X_s\>_H - 2\int_s^t\ _{V^*}\<Y_r, X_t\>_V\d r
-\|X_t\|_H^2-\|X_s\|_H^2 + 2 \<X_t, X_s\>_H \\
&= \|X_t\|_H^2 -\|X_s\|_H^2 -2 \int_s^t\ _{V^*}\<Y_r, X_t\>_V\d
r.\end{split}\end{equation*}

(b) As in \cite[Lemma 4.3]{KR}, 
we have
\beq\label{Tr10} \E\sup_{t\in
[0,T]}\|X_t\|_H^2<\infty.\end{equation} Indeed, by (\ref{Tr1}),
for any $t=t_i^l\in I_l\setminus \{0,T\}$ given in Lemma \ref{A1},

\begin{equation}\label{Tr11}
  \begin{split}
    &\|X_t\|_H^2 -\|X_0\|_H^2
    = \sum_{j=0}^{i-1} (\|X_{t_{j+1}^l}\|_H^2 -\|X_{t_j^l}\|_H^2 )\\
    &= 2\int_0^t\ _{V^*}\<Y_s, \tt X_s^l\>_V\d s
       + 2 \int_0^t \<\bar X_s^l, Z_s\d W_s\>_H
       + 2 \<X_0,\int_0^{t_1^l}Z_s \d W_s\>_H\\
    &\qquad + \sum_{j=0}^{i-1}
                 \big(
                   \|M_{t_{j+1}^l}-M_{t_j^l}\|_H^2
                   - \|X_{t_{j+1}^l}- X_{t_j^l}-M_{t_{j+1}^l}+M_{t_j^l}\|_H^2
                 \big).
  \end{split}
\end{equation}
By Remark~2.1(5), Lemma~\ref{A1} and {\bf (K)}(ii),

\beq\label{Tr12} \E \int_0^T |\ _{V^*}\<Y_s, \tt X_s^l\>_V|\d s\le
c \|\tt X^l\|_K\le c_1\end{equation} for some constant $c_1>0$
independent of $l.$ Moreover, by the Burkholder-Davis inequality,

\beq\label{Tr13}\beg{split} \E \sup_{t\le T} \bigg|\int_0^t \<\bar
X_s^l, Z_s\d W_s\>_H\bigg|&\le 3 \E \bigg(\int_0^T  \|\bar
X_s^l\|_H^2\d \<M\>_s\bigg)^{1/2}\\
&\le \ff 1 4 \E\sup_{k_l-1\ge j\ge 0} \|X_{t_j^l}\|_H^2 + 9 \E
\<M\>_T,\end{split}\end{equation}  where $\<M\>_t$ is the
increasing process part of $\|M_t\|_H^2.$ Finally, for all $l\geq 1$

\beq\label{Tr14}\limsup_{l\to\infty} \E \sum_{j= 0}^{k_l-1}
\|M_{t_{j+1}^l} -M_{t_j^l}\|_H^2 \le \E \<M\>_T.\end{equation}
Combining (\ref{Tr11})--(\ref{Tr14}), we obtain

$$\E \sup_{t\in I_l\setminus \{T\}} \|X_t\|_H^2 \le c_2$$
for some constant $c_2>0$ independent of $l$. Therefore, letting
$l\uparrow\infty$ and setting $I:=\cup_{l\ge 1}I_l\setminus\{T\},$
with $I_l$ as in Lemma \ref{A1}, we obtain

$$\E\sup_{t\in I} \|X_t\|_H^2<\infty.$$
Since for all $t\in [0,T]$

$$\sum_{j=1}^N \, _{V^*}\<X_t,e_j\>_V^2\uparrow \|X_t\|_H^2\ \text{as}\ N\uparrow \infty,$$
where as usual for $x\in V^*\setminus H$ we set $\|x\|_H:=\infty$,
it follows that $t\mapsto \|X_t\|_H$ is lower semicontinuous
$P$-a.s. Since $I$ is dense in $[0,T]$, we arrive at $\sup_{t\le
T}\|X_t\|_H^2 =\sup_{t\in I} \|X_t\|_H^2.$ Thus, (\ref{Tr10})
holds.

(c) Next, since $\sup_{t\in [0,T]}\|X_t\|_H<\infty$, the proof of
\cite[Lemma 4.5]{KR} implies

\beq\label{Tr2} \lim_{l\to\infty} \sup_{t\le
T}\bigg|\int_0^t\<X_s-\bar X_s^l,Z_s\d W_s\>_H\bigg|=0\ \
\text{in\ probability.}\end{equation}
We repeat the proof here for
completeness. We first note that because of (b) and its continuity
in $V^*$ the process $X$ is weakly continuous in $H$, and therefore,
since $\scr B(H)$ is generated by $H^*$, progressively measurable
as an $H$-valued process.
Hence, for
any $n\ge 1$ the process $P_n X_s$,
where $P_n$ is as defined in \eqref{Tr0},
is continuous in $H$ so that

$$\lim_{l\to\infty} \int_0^T \|P_n (X_s-\bar X_s^l)\|_H^2\d
\<M\>_s=0.$$ Therefore, it suffices to show that for any $\vv>0,$

\beg{equation}\label{Tr2'}\beg{split} &\lim_{n\to\infty}
\sup_{l\ge 1} P\bigg(\sup_{t\le T}\bigg|\int_0^t\<(1-P_n)\bar
X_s^l,Z_s\d
W_s\>_H\bigg|>\vv\bigg)=0,\\
& \lim_{n\to\infty} P\bigg(\sup_{t\le T}\bigg|\int_0^t\<(1-P_n)
X_s,Z_s\d W_s\>_H\bigg|>\vv\bigg)=0\;.\end{split}\end{equation}
For any $\dd\in (0,1)$ and $N>1$ by the Burkholder-Davis
inequality

\begin{equation*}\begin{split}
& P \bigg( \sup_{t\le T} \bigg| \int_0^t\<(1-P_n)\bar X_s^l,Z_s\d W_s\>_H \bigg| > \vv \bigg)
\le \ff{3\dd}{\vv} +P\bigg(\int_0^T\|\bar
X_s^l\|_H^2\d \<(1-P_n)M\>_s >\dd\bigg)\\
&\le \ff {3\dd}{\vv} + P\Big(\sup_{t\le T}\|X_t\|_H>N\Big) +\ff
{N^2}\dd \E\< (1-P_n)M\>_T.
\end{split}\end{equation*}
By letting
first $n\to \infty$, then $N\to\infty$ and finally $\dd\to 0$, we
prove the first equality in (\ref{Tr2'}). Similarly, the second
equality also holds.

(d) As in \cite[Lemma I.4.6]{KR}, we first prove (\ref{Ito}) for
$t\in I.$ So, fix $t\in I$. We may assume that $t\ne 0.$ In this
case for each sufficiently large $l\ge 1$ there exists a unique
$0<i<k_l$ such that $t= t_i^l.$ We have $X_{t_j^l}\in V$ a.s. for
all $j$.  By Lemma~\ref{A1}, Remark~2.1(5) and (\ref{Tr2}) the sum of
the first three terms in the right-hand side of (\ref{Tr11})
converges in probability to $2\int_0^t\ _{V^*}\<Y_s, \bar
X_s\>_V\d s+2\int_0^t \<X_s, Z_s\d W_s\>_H.$ Then

$$\|X_t\|_H^2 -\|X_0\|_H^2 = 2\int_0^t\ _{V^*}\<Y_s, \bar X_s\>_V\d s +
2 \int_0^t \<X_s, Z_s\d W_s\>_H +\<M\>_t-\vv_0,$$ where

$$\vv_0:= \lim_{l\to \infty}\sum_{t_{j+1}^l\le t}\|X_{t_{j+1}^l}-
X_{t_j^l}-M_{t_{j+1}^l}+ M_{t_j^l}\|_H^2$$ exists. So, to prove
(\ref{Ito}) for $t$ as above, it suffices to show that $\vv_0=0.$
Since for any $\varphi\in V$,

$$\<X_{t_{j+1}^l}-X_{t_j^l}-M_{t_{j+1}^l}+M_{t_j^l}, \varphi\>_H=
\int_{t_j^l}^{t_{j+1}^l}\ _{V^*}\<Y_s, \varphi\>_V \d s,$$ letting
$\tt M^l$ and $\bar M^l$ be defined as $\tt X^l$ and $\bar X^l$
respectively, for $M$ replacing  $X$, we obtain for every $n\ge 1$
\begin{equation*}\begin{split}
\vv_0=
\lim_{l\to\infty}\biggl(
  &\int_0^t\ _{V^*} \< Y_s,
                       \tt X_s^l - \bar X_s^l
                       - P_n (\tt M_s^l-\bar M_s^l)
                    \>_V\d s\\
  & -\bigbracket{X(t_1^l)-X(0)-M(t_1^l)+M(0), P_nM(0)-X(0)}_H\\
  & -\sum_{t_{j+1}^l\le t}
      \< X_{t_{j+1}^l}-X_{t_j^l}-M_{t_{j+1}^l}+M_{t_j^l},
         (1-P_n)(M_{t_{j+1}^l}-M_{t_j^l})\>_H
\biggr)\;.
\end{split}\end{equation*}
By the weak continuity of $X$ in $H$ the second term
converges to zero as $l\to\infty$.
Lemma~\ref{A1} and Remark~2.1(5)
imply that $\int_0^t \ _{V^*}\<Y_s, \tt X_s^l-\bar
X_s^l\>_V\d s \to 0$ in probability. Moreover, since $P_n M_s$ is
a continuous process in $V$, $\int_0^t \ _{V^*}\<Y_s, P_n(\tt
M_s^l-\bar M_s^l)\>_V\d s \to 0$ as $l\to\infty$. Thus,

\beg{equation*}\beg{split}\vv_0\le
&\lim_{l\to\infty}\Big(\sum_{t_{j+1}^l\le t}\|X_{t_{j+1}^l}-
X_{t_j^l}-M_{t_{j+1}^l}+
M_{t_j^l}\|_H^2\Big)^{1/2}\Big(\sum_{t_{j+1}^l\le t}\|(1-P_n)
(M_{t_{j+1}^l}-M_{t_j^l})\|_H^2\Big)^{1/2}\\= &\vv_0^{1/2}
\<(1-P_n)M\>_t^{1/2}, \end{split}\end{equation*} which goes to
zero as $n\to\infty$ since $M_t$ is a square-integrable martingale
in $H$. Therefore, $\vv_0=0.$

(e) Now, take $\OO'\in\scr F$ with full probability such that the
limit in (\ref{Tr2}) is a pointwise limit in $\OO'$ for some
subsequence (denoted again by $l\to\infty$)  and (\ref{Ito}) holds
for all $t\in I$  on $\OO'$. If $t\notin I$, for any $l\ge 1$
there exists a unique $j(l)<k_l$ such that $t\in (t_{j(l)}^l,
t_{j(l)+1}^l].$ Letting $t(l):= t_{j(l)}^l,$ we have $t(l)\uparrow
t$ as $l\uparrow \infty.$  By (\ref{Ito}) for $t\in I$, for any
$l>m$ we have, on $\OO'$ (since the above applies to $X-X_{t(m)}$
replacing $X$)

\beg{equation}\label{TR0}\beg{split} \|X_{t(l)}-X_{t(m)}\|_H^2 & =
2
\int_{t(m)}^{t(l)}\ _{V^*}\<Y_s, \bar X_s-X_{t(m)}\>_V\d s \\
&\qquad+ 2 \int_{t(m)}^{t(l)} \< X_s-X_{t(m)}, Z_s\d W_s\>_H +
\<M\>_{t(l)}
-\<M\>_{t(m)}\\
&= 2 \int_0^T 1_{[t(m), t(l)]}(s) _{V^*}\<Y_s, \bar X_s-\bar
X_s^m\>_V\d s \\
&\qquad + 2 \int_0^T1_{[t(m), t(l)]}(s) \< X_s-\bar X_s^m, Z_s\d
W_s\>_H+\<M\>_{t(l)} -\<M\>_{t(m)}.\end{split}\end{equation} Thus,
by the continuity of $\<M\>_t$ and (\ref{Tr2}) (holding pointwise
on $\OO'$), we have that

\beg{equation}\label{TR1} \lim_{m\to\infty}\sup_{l>m}\bigg\{ 2
\bigg|\int_0^T1_{[t(m), t(l)]}(s) \< \ X_s-\bar X_s^m, Z_s\d
W_s\>_H\bigg|+|\<M\>_{t(l)} -\<M\>_{t(m)}|\bigg\}=0\end{equation}
holding on $\OO'$.
Furthermore, by Lemma~\ref{A1} and by Remark~2.1(5),
selecting another
subsequence if necessary we have for some $\OO''\in \scr F$ with
full probability and $\OO''\subset \OO',$ that on $\OO''$

$$\lim_{m\to\infty}\int_0^T |\ _{V^*}\<Y_s, \bar X_s-\bar
X_s^m\>_V|\d s=0.$$ Since for all $t\notin I$

$$\sup_{l>m}\int_{t(m)}^{t(l)}|\ _{V^*}\<Y_s, \bar X_s-\bar
X_s^m\>_V|\d s\le \int_0^T|\ _{V^*}\<Y_s, \bar X_s-\bar X_s^m\>_V|\d
s,$$ we have that

$$\lim_{m\to\infty}\sup_{l>m}\int_{t(m)}^{t(l)} \ _{V^*}\<Y_s, \Bar X_s-\bar X_s^m\>_V\d
s=0$$ holds on $\OO''.$

 Combining this with
(\ref{TR0}) and (\ref{TR1}), we conclude that

$$\lim_{m\to\infty} \sup_{l\ge m}\|X_{t(l)}-X_{t(m)}\|_H^2=0$$
holds on $\OO''.$  Thus, $(X_{t(l)})_{l\in\mathbb N}$ converges in
$H$ on $\OO''$. Since we know that $X_{t(l)}\to X_t$ in $V^*$, it
converges to $X_t$ strongly in $H$ on $\OO''$. Therefore, by the
formula for $t(l)$ and letting $l\to\infty,$ we obtain (\ref{Ito})
on $\OO''$ also for all $t\notin I.$

Finally, since the right hand side of  (\ref{Ito}) is on $\OO''$
continuous in $t\in [0,T],$ so must be its left hand side, i.e.
$t\mapsto \|X_t\|_H$ is continuous on $[0,T].$ Therefore, the weak
continuity of $X_t$  in $H$ implies its strong continuity in $H$.
\end{proof}

\paragraph{Acknowledgement.} Jiagang Ren and Feng-Yu Wang (the corresponding author) would like to
thank the University of Bielefeld for a
very pleasant stay when this work was initiated.

\beg{thebibliography}{99}

\bibitem{Aronson} D.G. Aronson, \emph{The porous medium equation,} Lecture Notes Math.
Vol. 1224, Springer, Berlin, 1--46, 1986.

\bibitem{AP} D.G. Aronson and L.A. Peletier, \emph{Large time behaviour of solutions
of the porous medium equation in bounded domains,} J. Diff. Equ.
39(1981), 378--412.

\bibitem{BBDR} V. Barbu, V.I. Bogachev,  G. Da Prato and M. R\"ockner, \emph{Weak
solution to the stochastic porous medium  equations: the
degenerate case,} preprint.

\bibitem{BDR} V.I. Bogachev,  G. Da Prato and M. R\"ockner, \emph{Invariant measures
of   stochastic generalized porous medium equations,} to appear in
Dokl. Math.

\bibitem{DR1} G. Da Prato and M. R\"ockner, \emph{Weak solutions to stochastic
porous media equations,} J. Evolution Equ. 4(2004), 249--271.

\bibitem{DR2} G. Da Prato and M. R\"ockner, \emph{Invariant measures for a
stochastic porous medium equation,} preprint SNS, 2003; to appear
in Proceedings of Conference in Honour of K. It\^o, Kyoto, 2002.

\bibitem{DRRW} G. Da Prato, M. R\"ockner, B. L. Rozovskii and F.-Y.
Wang, \emph{Strong solutions to stochastic generalized porous
media equations: existence, uniqueness and ergodicity,} to appear
in Comm. Part. Diff. Equat.

\bibitem{DaZa} G. Da Prato and J. Zabczyk, \emph{ Stochastic Equations in
Infinite Dimensions,} Encyclopedia of Mathematics and its
Applications, Cambridge University Press. 1992.

\bibitem{DZ} A. Dembo and O. Zeitouni, \emph{ Large deviations
Techniques and Applications. Second Edition,} Springer, New York.
1998.

\bibitem{Doob} J. L. Doob, \emph{Stochastic Processes,} New York,
John Widely, 1953.

\bibitem{E} C. Ebmeyer, \emph{Regularity in Sobolev spaces for the
fast diffusion and the porous medium equation,} J. Math. Anal.
Appl. 307(2005), 134--152.

\bibitem{FOT94} M. Fukushima,  Y. Oshima and M. Takeda,
\emph{Dirichlet Forms and Symmetric Markov Processes,} Walter de
Gruyter, 1994.

\bibitem{Kim} J. U. Kim, \emph{On the stochastic porous medium
equation,}  J. Diff. Equat. 220(2006), 163--194.

\bibitem{K} N.V. Krylov, \emph{On Kolmogorov's equations for finite dimensional
diffusions,} Lecture Notes Math. 1715, 1--63, Berlin: Springer,
1999.

\bibitem{KR} N.V. Krylov and B.L. Rozovskii, \emph{Stochastic evolution equations,}
Translated from Itogi Naukii Tekhniki, Seriya Sovremennye Problemy
Matematiki 14(1979), 71--146, Plenum Publishing Corp. 1981.

\bibitem{P1} E. Pardoux, \emph{Sur des equations aux d\'eriv\'ees
partielles stochastiques monotones,} C. R. Acad. Sci. 275(1972),
A101--A103.

\bibitem{P2} E. Pardoux, \emph{Equations aux d\'eriv\'ees
partielles stochastiques non lineaires monotones: Etude de solutions
fortes de type Ito,} Th\'ese Doct. Sci. Math. Univ. Paris Sud. 1975.

\bibitem{RZ} M. M. Rao and Z. D. Ren, \emph{Applications  of Orlicz
Spaces,} New York: Marcel Dekker, 2002.

\bibitem{RWW} M. R\"ockner, F.-Y. Wang and L. Wu, Large deviations
for stochastic generalized porous media equations, preprint.

\bibitem{R} B. L. Rozovskii, \emph{Stochastic Evolution Systems:
Linear Theory and Applications to Nonlinear Filtering,} Dordrecht:
Kluwer Academic,  1990.

\bibitem{Zei90} E. Zeidler, \emph{Nonlinear Functional Analysis
and Its Applications II/B,}  New York: Springer, 1990.

\end{thebibliography}

\end{document}